\documentclass{amsart}
\usepackage{amsthm}
\usepackage{amsmath,amssymb,amsfonts}

\usepackage{xspace,enumerate}
\usepackage{burdges}
\usepackage{calc}

\newcommand\Uir{U_{0,r}}

\newcommand\elabpgrps{\mathcal{E}}
\newcommand\components{\mathcal{E}}

\newcommand\reB{\tilde{B}}
\newcommand\reU{\tilde{U}}
\newcommand\reE{\tilde{E}}
\newcommand\reF{\tilde{F}}
\newcommand\reGamma{\tilde{\Gamma}}
\newcommand\recomponents{\tilde{\components}}

\begin{document}

\title[Signalizers and balance]%
{Signalizers and balance in groups \\ of finite Morley rank}
\author{Jeffrey Burdges}
\thanks{Burdges was supported by NSF grant DMS-0100794, and
 Deutsche Forschungsgemeinschaft grant Te 242/3-1.}
\address{Jeffrey Burdges\\
School of Mathematics\\
The University of Manchester\\
PO Box 88, Sackville St.\\
Manchester M60 1QD, England}
\email{Jeffrey.Burdges@manchester.ac.uk}


\begin{abstract}
There is a longstanding conjecture, due to Gregory Cherlin and Boris Zilber,
that all simple groups of finite Morley rank are simple algebraic groups.
The most successful approach to this conjecture has been Borovik's program
analyzing a minimal counterexample, or simple $K^*$-group.
We show that a simple $K^*$-group of finite {M}orley rank and odd type is
either aalgebraic of else has Prufer rank at most two.
This result signifies a switch from the general methods used to handle
large groups, to the specilized methods which must be used to identify
 $\PSL_2$, $\PSL_3$, $\PSp_4$, and $G_2$.
\end{abstract}

\maketitle

The Algebraicity Conjecture for simple groups of finite Morley rank,
also known as the Cherlin-Zilber conjecture, states that simple groups
of finite Morley rank are simple algebraic groups over algebraically
closed fields.  In the last 15 years, the main line of attack on this problem
has been the Borovik program of transferring methods from finite group
theory.  This program has led to considerable progress; however,
the conjecture itself remains decidedly open.
We divide groups of finite Morley rank into four types, odd, even, mixed,
and degenerate, according to the structure of their Sylow 2-subgroups.
For even and mixed type the Algebraicity Conjecture has been proven,
and connected degenerate type groups are now known to have trivial
Sylow 2-subgroups \cite{BBC}.
The present paper is part of the program to analyze a minimal
counterexample to the conjecture in {\em odd type}, where the
Sylow 2-subgroup is divisible-abelian-by-finite.
It is the final paper in a sequence proving that such a minimal
counterexample, or simple nonalgebraic $K^*$-group, has
\Prufer 2-rank at most two.

\begin{namedtheorem}{High \Prufer Rank Theorem}
A simple $K^*$-group of finite Morley rank with \Prufer 2-rank at least three
is algebraic.
\end{namedtheorem}

This will be a consequence of the following so-called trichotomy,
 which is proved in the present paper.
Here the traditional term ``trichotomy'' refers to the fact that there is also
the \Prufer 2-rank $\leq 2$ case, which is largely unexplored at present.

\begin{namedtheorem}{Generic Trichotomy Theorem}\label{generictrico_intro}
Let $G$ be a simple $K^*$-group of finite Morley rank and odd type
with \Prufer 2-rank $\geq 3$.  Then either
\begin{conclusions}
\item $G$ has a proper 2-generated core, or
\item $G$ is an algebraic group over an algebraically closed field of characteristic \noteql2.
\end{conclusions}
\end{namedtheorem}

\noindent The High \Prufer Rank Theorem then follows
by applying the next two results.
% (and recalling that \Prufer rank is less then normal rank).

\begin{namedtheorem}{Strong Embedding Theorem}[{\cite{BBN04}}]
Let $G$ be a simple $K^*$-group of finite Morley rank and odd type,
 with normal 2-rank $\geq3$ and \Prufer 2-rank $\geq2$.
% ---both of which follow from \Prufer 2-rank $\geq 3$.
Suppose that $G$ has a proper 2-generated core $M$. % $M = \Gamma_{S,2}(G)$.
Then $G$ is a minimal connected simple group,
 and $M$ is strongly embedded.
\end{namedtheorem}

\begin{namedtheorem}{Minimal Simple Theorem}[{\cite{BCJ}}]
Let $G$ be a minimal connected simple group of finite Morley rank and of odd type.
Suppose that $G$ contains a proper definable strongly embedded subgroup $M$.
Then $G$ has \Prufer 2-rank one.
\end{namedtheorem}

It may seem odd that the first of these results is appearing last.
In fact, an earlier version of the trichotomy theorem began % initiated
this sequence of developments.
Namely, Borovik first proved the trichotomy theorem
under a {\em tameness} assumption in \cite{Bo95}, and
the present author had explored eliminating tameness in \cite{Bu03}.
In \cite{Bo95}, Borovik produces the proper 2-generated core
with a tame nilpotent signalizer functor theorem
 \cite[\qTheorem B.30]{BN} (see also \cite[\qTheorem 6.2]{Bu03}),
 an approach mirrored in the present paper.
In \cite{Bu03}, we show that the ``most unipotent part'' of a
solvable signalizer functor is a nilpotent signalizer functor.
This was believed to quickly eliminate tameness from \cite{Bo95}.
However, more careful investigations revealed that obtaining
a signalizer functor remained problematic.

In \S\ref{subsec:rebalancing} of the present paper,
we resolve this difficulty by constructing signalizer functors of
a ``sufficiently unipotent'' reduced rank.  The most serious obstacle
is explained in Example \ref{rebalancing_example}.
Our approach forces subsequent analysis to restrict itself to
components of the centralizers of involutions which involve
sufficiently large fields, a worrying but ultimately harmless
restriction.  Indeed, all complexities introduced by this approach 
are dispensed with in \S\ref{subsec:rebalancing}.
This seems to be a different approach from that used by finite
group theorists, who work with so-called {\em weakly balanced}
signalizer functors \cite[\S29]{GLS1};
a similar method might work here as well.
We call our approach partial balance.

\smallskip

The first section of this article covers necessary background material,
 including the definitions of a signalizer functor and the 2-generated core.

The second section contains the delicate definitions of partial balance,
and of the associated family $\recomponents_X$ of components
from the centralizers of involutions.  This section also contains a
version of Asar's theorem (Theorem \ref{rebalanced_Asar1})
 which states that $\recomponents_X \neq \emptyset$,
as well as a criterion for $\gen{\recomponents_X} = G$
 (Theorem \ref{rebalanced_Asar2}).
Borovik's earlier unpublished work on the analysis of Lie rank two
 components \cite{Bo03} has heavily influenced this final result,
although partial balance has given these results a more technical flavor.

The third section provides a suitable version of Berkman and Borovik's
 Generic Identification Theorem \cite{BB04}.
It is the role of section two to verify the two hypotheses of this argument:
 reductivity for, and generation by, the centralizers of involutions. 
Our partial balance approach provides only a weak form of
 the reductivity hypothesis, which necessitates some alterations in
 the proof of the Generic Identification Theorem.
All such critical changes are confined to \S\ref{subsec:RootGroups},
but there are important modifications throughout \S\ref{sec:generictrico}.
The reader unfamiliar with the Generic Identification Theorem
should consider exploring \S\ref{sec:generictrico} before
 \S\ref{sec:preliminaries} or \S\ref{sec:reBnC}.

\smallskip

This is by no means the end of the story.
The \Prufer 2-rank $\geq 3$ hypothesis used here is weaker than the
normal 2-rank $\geq 3$ hypothesis originally used by Borovik \cite{Bo95}.
As part of the ongoing program in odd type, \cite{BB07} will show
 that Borovik's original trichotomy holds, without the tameness hypothesis.
We view \cite{BB07} as a bridge between the ``generic case'' which is
 treated here, and the ``quasi-thin'' case
 (the identification of $\PSp_4$, $\G_2$, and $\PSL_3$).

\section{Preliminaries}\label{sec:preliminaries}

This first section recalls various definitions and facts which are used
throughout \S\ref{sec:reBnC}, and less pervasively in \S\ref{sec:generictrico}.
% Our attention necessarily focuses on the 2-generated core and signalizer functors.

\subsection{$K$-groups}

We proceed, in this paper, by analyzing a so-called simple $K^*$-group
of finite Morley rank.  A $K^*$-group is a group whose {\em proper}
definable simple sections are all algebraic.  Similarly, a $K$-group is a
group whose definable simple sections are all algebraic.
So the proper subgroups of our $K^*$-group are clearly $K$-groups.
One major $K$-group fact used throughout this article is the following
generation principle.

\begin{fact}[{\cite[\qTheorem 5.14]{Bo95}; see also \cite[\qTheorem 3.25]{BuPhd}}]
\label{CgenerationKU2}
Let $G$ be a connected $K$-group of finite Morley rank and odd type.
Let $V$ be a four-subgroup acting definably on $G$.  Then
$$ G = \gen{ C^\o_G(v) \mid v\in V^\# } $$
\end{fact} 

% The proof of this fact began with solvable group theory, as well as
% the facts of the next subsection, and ended by looking at the Shur multipliers
% and subsystem subgroups of a simple algebraic group.

Another major $K$-group fact worth recalling at the outset is the following
``reductivity'' criterion, which requires two definitions.

\begin{definition}
A quasisimple subnormal subgroup of a group $G$ is called a
{\em component} of $G$ (see \cite[p.~118 (2)]{BN}).
We define $E(G)$ to be
 the connected part of the product of components of $G$,
 or equivalently the product of the components of $G^\o$
(see \cite[\qLemma 7.10iv]{BN}).
Such components are normal in $G^\o$ by \cite[\qLemma 7.1iii]{BN},
 and indeed $E(G) \normal G$.
\end{definition}

\begin{definition}
The {\em odd part} $O(G)$ of a group $G$ of finite Morley rank is
 the maximal definable connected normal $2^\perp$-subgroup of $G$.
\end{definition}

Clearly $O(G)$ is solvable if $G$ is a $K$-group.

\begin{fact}[{\cite[\qTheorem 5.12]{Bo95}}]\label{KgrpO}
Let $H$ be a connected $K$-group of finite Morley rank and odd type
with $O(H) = 1$.  Then $H = F^\o(H) * E(H)$ is isomorphic to
a central product of quasisimple algebraic groups over
algebraically closed fields of characteristic \noteql2 and
of a definable normal divisible abelian group $F^\o(H)$.
\end{fact}

This fact motivates signalizer functor theory, whose goal is to show that
$O(H) = 1$, or something similar, when $H$ is the centralizer of an involution.
The version we aim at here is Corollary \ref{KgrpO_reU} below.

\subsection{Algebraic groups}

A key tool in our program is the fact that a group of finite Morley rank
 acting faithfully as a group of automorphisms of an algebraic group
 must itself be algebraic.

\begin{definition}
Given an algebraic group $G$, a maximal torus $T$ of $G$, and a
Borel subgroup $B$ of $G$ which contains $T$, we define the group
$\Gamma$ of {\em graph automorphisms} associated to $T$ and $B$,
to be the group of algebraic automorphisms of $G$ which normalize
both $T$ and $B$.
\end{definition}

\begin{fact}[{\cite[\qTheorem 8.4]{BN}}]\label{autalg}
Let $G \rtimes H$ be a group of finite Morley rank where $G$ and
$H$ are definable, $G$ an infinite quasi-simple algebraic group over an
algebraically closed field, and $C_H(G)$ is trivial.  Then, viewing $H$
as a subgroup of $\Aut(G)$, we have $H\leq \Inn(G)\Gamma$, where
$\Inn(G)$ is the group of inner automorphisms of $G$ and $\Gamma$
is the group of graph automorphisms of $G$, relative to a fixed choice
of Borel subgroup $B$ and maximal torus $T$ contained in $B$.
\end{fact}

An algebraic group is said to be {\em reductive} if it has no unipotent
radical.  Such a group is a central product of semisimple algebraic
groups and algebraic tori.  The centralizer of an involution in a
reductive algebraic group over a field of characteristic $\neq 2$
is itself reductive.

% \begin{fact}[{\cite[\qTheorem 3.5.4]{Carter93}}]\label{Creductive}
% Let $G$ be a connected reductive algebraic group.  Then the
% group $C^\o_G(s)$ is reductive for any semisimple automorphism $s$ of $G$.
% \end{fact}

\begin{fact}[{\cite[Theorem 8.1]{St2}}]\label{Creductive2}
Let $G$ be a quasisimple algebraic group over an algebraically closed field.
Let $\phi$ be an algebraic automorphism of $G$ whose order is finite and
relatively prime to the characteristic of the field.
Then $C^\o_G(\phi)$ is nontrivial and reductive.
\end{fact}

% \begin{fact}\label{Cquotient_fin}
% Let $G = H T$ be a group of finite Morley rank.  Suppose that
%  $Q \normaleq H$ is a finite $T$-invariant subgroup which is
% central in $H$.  Then $$ C^\o_{H/Q}(T) = C^\o_H(T)Q/Q
%   \quad\textrm{and}\quad  [C_{H/Q}(t) : C^\o_{H/Q}(t)] \mid \abs{Q}
%   \quad\textrm{for any $t\in T$}\mathperiod $$
% \end{fact}  %% Needs H normal in G

% We may assume that $G$ is its own universal central extension
%  since $C^\o_{G/Z(G)}(\phi) = C^\o_G(T)/Z(G)$
% by Fact \ref{Cquotient_fin}.  % was by \cite[\qFact 3.2]{Bu03}.
\begin{proof}
We know both that $C_H(\phi) \leq C_H(\phi \bmod Q)$
 as well as $[\phi,C_H(\phi \bmod Q)] \leq Q$.
There is a homomorphism $C_H(\phi \bmod Q) \to Q$ given by $x \mapsto [\phi,x]$.
As $Z(G)$ is finite,
 $t$ centralizes $C^\o_H(\phi \bmod Q)$, as desired.
So $C^\o_{G/Z(G)}(\phi) = C^\o_G(T)/Z(G)$
We may therefore assume that $G$ is its own universal central extension.

Since $\phi$ is algebraic and has finite order, $G \rtimes \gen{\phi}$ is
an algebraic group which contains $\phi$ as an inner automorphism.
Since the order of $\phi$ is finite and relatively prime to the characteristic,
$\phi$ is a semisimple automorphism of $G$.
So the result follows from Theorem 8.1 of \cite{St2}.
\end{proof}

% \begin{proof}
% We may assume that $G$ is its own universal central extension
%  since $$ C^\o_{G/Z(G)}(\phi) = C^\o_G(\phi)/Z(G)\mathperiod $$
% % by Fact \ref{Cquotient_fin}.
% Since $\phi$ is algebraic and has finite order,
%  $G \rtimes \gen{\phi}$ is an algebraic group which contains
%  $\phi$ as an inner automorphism.
% Since the order of $\phi$ is finite and relatively prime to the characteristic,
%  $\phi$ is a semisimple automorphism of $G$.
% So the result follows from Theorem 8.1 of \cite{St2}.
% \end{proof}

More specialized facts about algebraic groups will appear in \S\ref{sec:generictrico}.

\subsection{Unipotent groups}

While there is no intrinsic definition of unipotence in a group of 
finite Morley rank, there are various analogs of the ``unipotent radical'':
the Fitting subgroup, the $p$-unipotent operators $U_p$, for $p$ prime,
and their ``characteristic zero'' analogs $\Uir$ from \cite{Bu03,BuPhd}.
We recall their definitions below.

\begin{definition}
The {\em Fitting subgroup} $F(G)$ of a group $G$
 of finite Morley rank is the subgroup generated
 by all its nilpotent normal subgroups.
\end{definition}

The Fitting subgroup is itself nilpotent and definable
 \cite[\qTheorem 7.3]{Bel87,Ne91,BN},
and serves as a rough notion of unipotence in some contexts.
However, the Fitting subgroup of a solvable group $H$
may not be contained in the Fitting subgroup of a solvable group
containing $H$. % so ?it? is not a robust notion.

\begin{definition}
A connected definable $p$-subgroup of bounded exponent in
a group $H$ of finite Morley rank is said to be {\it $p$-unipotent}.
We write $U_p(H)$ for the subgroup generated by
 all $p$-unipotent subgroups of $H$.
\end{definition}

Clearly $U_p(H)$ need not be solvable when $H$ is a non-solvable 
algebraic group in characteristic $p$; however, a $p$-unipotent
$K$-group is solvable, and hence nilpotent by the following.

\begin{fact}
[{\cite[\qCorollary 2.16]{CJ01}; \cite[\qFact 2.36]{ABC97}}]
\label{Upnilpotence} % (Existence and nilpotence of the $p$-unipotent radical)
Let $H$ be a connected solvable group of finite Morley rank.
Then $U_p(H) \leq F^\o(H)$ is itself $p$-unipotent, and hence nilpotent.
\end{fact}

Thus the $p$-unipotent radical $U_p$ will automatically behave
well, inside a solvable group.  Its only weakness is that it may be trivial.

\begin{fact}[{\cite[\qTheorem 9.29 \& \S6.4]{BN}}]\label{Sylow_conplus}
Let $G$ be a connected solvable group of finite Morley rank. 
Then a Sylow $p$-subgroup $P$ of $G$ is connected, and
 $P = U_p(G) * T$ for a divisible abelian $p$-group $T$.
\end{fact}

The present paper relies on the theory of ``characteristic zero''
unipotence introduced in \cite{Bu03}.
We now turn our attention to this definition,
 as well as some facts from \cite{Bu03,Bu05a,BuPhd}.

\begin{definition}
We say that a connected abelian group of finite Morley rank is
{\it indecomposable} if it has a unique maximal proper
definable connected subgroup, denoted $J(A)$ (see \cite[\qLemma 2.4]{Bu03}).
We define the {\it reduced rank} $\rr(A)$ of
a definable indecomposable abelian group $A$
to be the Morley rank of the quotient $A/J(A)$,
i.e.\ $\rr(A) = \rk(A/J(A))$.
For a group $G$ of finite Morley rank, and any integer $r$, we define
$$ \Uir(G) = \Genst{A \leq G}{%
\parbox{\widthof{$A$ is a definable indecomposable group,}}%
{$A$ is a definable indecomposable group, \\
\hspace*{10pt} $\rr(A) = r$, and $A/J(A)$ is torsion-free}}\mathperiod $$
We say that $G$ is a {\it $\Uir$-group} (alternatively {\it $(0,r)$-unipotent} group)
if $U_{0,r}(G)=G$.  We also set $\rr_0(G) = \max \{r \mid \Uir(G) \neq 1 \}$.
\end{definition}

We view the reduced rank parameter $r$ as a {\em scale of unipotence},
 with larger values being more unipotent.
By the following fact, analogous to Fact \ref{Upnilpotence},
 the ``most unipotent'' groups, in this scale, are nilpotent.

\begin{fact}%
[{\cite[\qTheorem 2.21]{BuPhd}; \cite[\qTheorem 2.16]{Bu03}}]
\label{nilpotence}
Let $H$ be a connected solvable group of finite Morley rank.
Then $U_{0,\rr_0(H)}(H) \leq F(H)$.
\end{fact}

\begin{fact}[{\cite[\qCorollary 4.6]{Bu05a}}]\label{quicknilpotence2}
Let $G = H T$ be a group of finite Morley rank,
 with $H$ and $T$ definable and nilpotent, and $H \normal G$.
Suppose that $T$ is a $\Uir$-group for some $r \geq \rr_0(H)$.
Then $G$ is nilpotent.
\end{fact}

A {\it good torus} is a divisible abelian group of finite Morley rank 
whose definable connected subgroups are the definable hulls of their torsion.
We arrive at a good torus when all our various notions of unipotence are trivial.

\begin{fact}%
[{\cite[\qTheorem 2.19]{BuPhd}; \cite[\qTheorem 2.15]{Bu03}}]
\label{Udichotomy}
Let $H$ be a connected solvable group of finite Morley rank.  Suppose
$U_p(H) = 1$ for all $p$ prime, and $U_{0,\rr_0(H)}(H) = 1$.
Then $H$ is a good torus.
\end{fact}

In a similar vein, the notion of $(0,r)$-unipotence provides a
 useful decomposition of a nilpotent group.

\begin{fact}%
[{\cite[\qCorollary 3.6]{Bu05a}; \cite[\qTheorem 2.31]{BuPhd}}]
\label{nildecomp}
Let $G$ be a connected nilpotent group of finite Morley rank.
Then $G = D * B$ is a central product of definable characteristic
 subgroups $D,B \leq G$ where $D$ is divisible and
 $B$ is connected of bounded exponent.
Let $T$ be the torsion part of $D$.
Then we have decompositions of $D$ and $B$ as follows.
\begin{eqnarray*}
D &=& d(T) * U_{0,1}(G) * U_{0,2}(G) * \cdots \\
B &=& U_2(G) \times U_3(G) \times U_5(G) \times \cdots
\end{eqnarray*}
\end{fact}

The next fact tells us when $q$-unipotence is preserved by taking centralizers,
 a fact used to produce a signalizer functor in Lemma \ref{rebalanced_Uq} below.

\begin{fact}[{\cite[Fact 3.4]{Bu03}; \cite{ABCCdraft}}]\label{Cconnected}
Let $G$ be a connected solvable $p^\perp$-group of finite Morley rank,
and let $P$ be a finite $p$-group of definable automorphisms of $G$.
Then $C_G(P)$ is connected.
\end{fact}

There is also a ``characteristic zero'' analog of the foregoing.

\begin{fact}[{\cite[Lemma 3.6]{Bu03}}]\label{Cunipotent}
Let $G$ be a nilpotent $(0,r)$-unipotent $p^\perp$-group of finite Morley rank,
and let $P$ be finite $p$-group of definable automorphisms of $G$.
Then $C_G(P)$ is $(0,r)$-unipotent.
\end{fact}

In a similar vein, commutator subgroups of connected or $(0,r)$-unipotent
 groups tend to retain these properties.

\begin{fact}[{\cite[\qCorollary 5.29]{BN}}]\label{concommutator}
Let $H$ be a definable connected subgroup of a group $G$ of finite Morley rank
and let $X \subset G$ be any subset of $G$.  Then the group $[H,X]$ is definable
and connected.
\end{fact} 

\begin{fact}[{\cite[\qCorollary 3.6]{Bu05a}}]\label{Unilcommutator}
Let $G$ be a solvable group of finite Morley rank, let $S \subset G$ be
any subset, and let $H$ be a nilpotent $\Uir$-group which is normal in $G$.
Then $[H,S] \leq H$ is a $\Uir$-group. 
\end{fact}

\subsection{2-Local structure}

As the goal of our project is to constrain the 2-local structure,
we need a few parameters to measure the complexity of a Sylow 2-subgroup.
We define the {\em 2-rank} $m(G)$ of a group $G$ to be the
maximum rank of its elementary abelian 2-subgroups.
The {\em \Prufer 2-rank} $\pr(G)$ is the maximum $k$ such that
there is a \Prufer 2-subgroup $\Z(2^\infty)^k$ inside $G$, and
the {\em normal 2-rank} $n(G)$ is the maximum 2-rank of a normal
elementary abelian 2-subgroup of $G$.  In an odd type group
of finite Morley rank, these various ranks are all finite,
and we have $$ m(G) \geq n(G) \geq \pr(G)\mathperiod $$
These notions are well-defined because the Sylow 2-subgroups of
a group of finite Morley rank are conjugate \cite[\qTheorem 10.11]{BP,BN}.

We use $\elabpgrps_k(H)$ to denote the set of elementary abelian
2-subgroups $U \leq H$ with $m(U) \geq k$.
We give $\elabpgrps_2(H)$ a {\em graph structure} by placing
an edge between $U,V\in \elabpgrps_2(H)$ whenever $[U,V] = 1$.
We say $H$ is {\em 2-connected} if the graph $\elabpgrps_2(H)$
is connected, and we refer to the components of $\elabpgrps_2(H)$
as {\em 2-connected components} otherwise.

\begin{fact}[{compare \cite[46.2]{Asch}}]\label{Asch46.2}
Let $S$ be a locally finite 2-group.  Then
\begin{conclusions}
\item  If $m(S) > 2$ then the graph $\elabpgrps_2(S)$ has a
       unique nonsingleton 2-connected component given by
$$ \elabpgrps^0_2(S) := \{ X\in \elabpgrps_2(S) : m(C_S(X))>2 \}, \textrm{and} $$
       $\elabpgrps^0_2(S)$ contains any $X\in \elabpgrps_2(S)$
       with $X\normal S$.
\item  If $n(S) > 2$ then $S$ is 2-connected.
\end{conclusions}
\end{fact}

\begin{proof}
Since $S$ is locally finite, this reduces to the finite case,
 found in \cite[46.2]{Asch}.
\end{proof}

\subsection{Proper 2-generated core}

\begin{definition}
Consider a group $G$ of finite Morley rank
 and a 2-subgroup $S$ of $G$ with $m(S) \geq 3$.
We define the {\em 2-generated core} $\Gamma_{S,2}(G)$ of $G$
(associated to $S$) to be the definable hull of the group generated by
all normalizers of groups in $\elabpgrps_2(S)$:
$$ \Gamma_{S,2}(G) = d(\gen{ N_G(U) : U\in \elabpgrps_2(S) }) $$
We also define the {\em weak 2-generated core} $\Gamma^0_{S,2}(G)$
of $G$ (associated to $S$) to be the definable hull of all normalizers of
groups in the nonsingleton 2-connected component $\elabpgrps^0_2(S)$.
$$ \Gamma^0_{S,2}(G) =
   d(\gen{ N_G(U) : U\in \elabpgrps_2(S), m(C_S(U))>2 })\mathperiod $$
We say that $G$ has a {\em proper 2-generated core},
 or a {\em proper weak 2-generated core},
 when, for a Sylow 2-subgroup $S$,
 $\Gamma_{S,2}(G) < G$ or $\Gamma^0_{S,2}(G) < G$, respectively.
\end{definition}
Both notions of 2-generated core are well-defined,
 by the conjugacy of Sylow 2-subgroups.
By Fact \ref{Asch46.2}-2, the 2-generated core and the
weak 2-generated core coincide when $n(G) \geq 3$,
 as is the case for much of the rest of this paper.
When they differ, the weak 2-generated core is the more useful notion.

For an elementary abelian 2-group $V$ acting definably on $G$,
we define $\Gamma_V(G)$ to be the group generated by the
connected centralizers of involutions in $V$.
$$ \Gamma_V(H) = \gen{C^\o_H(v) : v\in V^\#}\mathperiod $$

\begin{proposition}\label{GITcondA}
Let $G$ be a simple $K^*$-group of finite Morley rank and odd type,
with $m(G)\geq3$, and let $S$ be a Sylow 2-subgroup of $G$.
Suppose that $\Gamma_E(G) < G$ for some $E \in \elabpgrps^0_2(S)$.
Then $G$ has a proper weak 2-generated core.
\end{proposition}

% For this, we need to know that $C^\o_G(i) \neq 1$ for $i$ an involution.
%
% \begin{fact}[{\cite[Ex. 13 \& 15 p.~79]{BN}}]\label{centinv}
% Let $G$ be a connected group of finite Morley rank and
% let $j : G \to G$ be a definable involutive automorphism.
% If the centralizer $C_G(j)$ is finite, then $G$ is abelian
%  and $j$ acts on $G$ by inversion.
% \end{fact}

This depends on a lemma.

\begin{lemma}\label{preGITcondA} % [{cf. \cite[Lemma 3.1]{Bo03}}]
Let $G$ be a simple $K^*$-group of finite Morley rank and odd type.
Then $\Gamma_U(G) = \Gamma_V(G)$ for any $U,V$ in the same connected
component of the graph $\elabpgrps_2(G)$.
\end{lemma}

\begin{proof}
It is enough to prove the result for $U,V$ with $[U,V] = 1$.
For any $v\in V^\#$, simplicity implies that $C^\o_G(v)$ is
a proper subgroup of $G$, and hence a $K$-group.
Since $U$ normalizes $C^\o_G(v)$,
 $C^\o_G(v) = \Gamma_U(C^\o_G(v))$ by Fact \ref{CgenerationKU2}.
So $\Gamma_V(G) \leq \Gamma_U(G)$, and the result follows by symmetry.
\end{proof}

\smallskip

\begin{proof}[of Proposition \ref{GITcondA}]
We may assume $E \leq S$ by conjugacy of Sylow 2-subgroups.
Since involutions of $G$ have infinite centralizers by \cite[Ex. 13 \& 15 p.~79]{BN}, % Fact \ref{centinv},
 the result will follow from the following claim, and simplicity.
$$ \Gamma^0_{S,2}(G) \leq N_G(\Gamma_E(G))
  \quad\textrm{for any $E\in \elabpgrps^0_2(S)$} $$
By Lemma \ref{preGITcondA} and Fact \ref{Asch46.2}-1,
 $\Gamma_E(G) = \Gamma_U(G)$ for any $U\in \elabpgrps^0_2(S)$.
For any $U\in \elabpgrps^0_2(S)$,
 $$ N_G(U) \leq N_G(\Gamma_U(G)) = N_G(\Gamma_E(G))\mathperiod $$
Thus $\Gamma^0_{S,2}(G) \leq N_G(\Gamma_E(G))$, as desired.
\end{proof}

We will encounter a variation of the preceding in the next section
 (see Lemma \ref{rebalanced_preGITcondA} and
 Proposition \ref{rebalanced_GITcondA}).

The following black hole principle for proper 2-generated cores reverses
the roles of the subgroups $\Gamma^0_{S,2}(G)$ and $\Gamma_E(G)$
 in Proposition \ref{GITcondA}.

\begin{lemma}\label{weak_CconM}
Let $G$ be an infinite simple $K^*$-group of finite Morley rank and odd type,
 and let $S$ be a 2-subgroup of $G$ satisfying $m(S)\geq3$.
Then $C^\o_G(x) \leq \Gamma^0_{S,2}(G)$ for every
 $x\in I(S)$ with $m(C_S(x)) > 2$. % $x \in \bigcup \elabpgrps^0_2(S) \setminus \{1\}$.
\end{lemma}

\begin{proof}
There is an $E\in \elabpgrps^0_2(S)$ with $x\in E$ and $m(E)\geq3$
by Fact \ref{Asch46.2}-1.  So there is an $E_1\in \elabpgrps^0_2(S)$
with $E_1 \leq E$ and $E_1 \cap \gen{x} = 1$.
For any $y\in E_1^\#$, we have $C_{C_G(x)}(y) \leq C_G(y,x)$ and
$\gen{y,x} \in \elabpgrps^0_2(S)$.
By simplicity, Fact \ref{CgenerationKU2} yields
\[ C^\o_G(x) = \Gamma_{E_1}(C^\o_G(x)) \leq \Gamma^0_{S,2}(G) \]
\end{proof}

In particular, given a simple $K^*$-group $G$,
 Lemma \ref{weak_CconM} says
$$ \Gamma_E(G) \leq \Gamma_{E,2}(G)
    \quad\textrm{for any $E\in \elabpgrps_3(G)$\mathperiod} $$
Now Proposition \ref{GITcondA} and Lemma \ref{weak_CconM}
 yield the following.
 
\begin{proposition}\label{twogencore_blowup}
Let $G$ be a simple $K^*$-group of finite Morley rank and odd type,
with $m(G)\geq3$, and let $S$ be a Sylow 2-subgroup of $G$.
If $\Gamma_{E,2}(G) < G$ for some $E \in \elabpgrps_3(G)$,
then $G$ has a proper weak 2-generated core, i.e.\ $\Gamma^0_{S,2}(G) < G$.
\end{proposition}

% \begin{proof}
% By Lemma \ref{weak_CconM}, $\Gamma_E(G) \leq \Gamma_{E,2}(G) < G$.
% By Proposition \ref{GITcondA},
%  $G$ has a proper weak 2-generated core $\Gamma^0_{S,2}(G) < G$.
%  \end{proof}

\subsection{Signalizer functors}

Signalizer functors are used in both the finite case and in the
finite Morley rank case to produce a dichotomy between
a proper 2-generated core, and a reductivity condition for
centralizers of involutions.

\begin{definition}
Consider a group $G$ of finite Morley rank, and an elementary abelian
2-subgroup $E\in \elabpgrps_3(G)$.
An {\em $E$-signalizer functor} on $G$ is a family $\{\theta(s)\}_{s\in E^\#}$
 of definable $E$-invariant $2^\perp$-subgroups of $G$ satisfying:
\begin{hypotheses}
\item  $\theta(s) \normaleq C_G(s)$ for each $s\in E^\#$.
\item  $\theta(s) \cap C_G(t) \leq \theta(t)$ for any $s,t \in E^\#$.
\end{hypotheses}
\end{definition}
% \cite[p.\ 124]{GLS1}

We observe that the second condition is equivalent to the ``balance'' condition
$$ \theta(s) \cap C_G(t) = \theta(t) \cap C_G(s)  \quad
 \textrm{for any $s,t \in E^\#$\mathperiod} $$
In practice, we will only be interested in signalizer functors satisfying
the following stronger invariance condition,
 which is used to produce a proper 2-generated core.
\begin{texteqn}{\tag{$\dag$}}
$\theta(s)^g = \theta(s^g)$ for all $s\in E^\#$ and
 all $g\in G$ for which $s^g\in E$.
\end{texteqn}
As one would expect, we say $\theta$ is a {\em connected} or
{\em nilpotent} signalizer functor if the groups $\theta(s)$ are
connected or nilpotent, respectively, for {\em all} $s\in E^\#$.

We now show that signalizer functors yield a proper weak 2-generated core.

\begin{theorem}\label{twogencore}
Let $G$ be a simple $K^*$-group of finite Morley rank and odd type,
 and let $S$ be a Sylow 2-subgroup of $G$.
Suppose that, for some $E\in \elabpgrps_3(S)$, $G$ admits
 a nontrivial connected nilpotent $E$-signalizer functor $\theta$ satisfying $(\dag)$
% \begin{texteqn}{\tag{$\star$}}
% $\theta(s)^g = \theta(s^g)$ for all $s\in E^\#$ and
%  all $g\in G$ for which $s^g\in E$.
% \end{texteqn}
Then $G$ has a proper weak 2-generated core.
\end{theorem}

The key fact underlying this result is the Nilpotent Signalizer Functor Theorem.

\begin{definition}
We say that an $E$-signalizer functor on a group $G$
of finite Morley rank is {\em complete} if:
\begin{hypotheses}
\item $\theta(E) = \gen{ \theta(s) : s \in E^\# }$ is a solvable $p^\perp$-group, and
\item $\theta(s) = C_{\theta(E)}(s)$ for any $s \in E^\#$.
\end{hypotheses}
\end{definition}

\begin{namedtheorem}{Nilpotent Signalizer Functor Theorem}
[{\cite[Thm.\ B.30]{Bo95,Bu03,BN}}]
\label{named:Ncompleteness}
Let $G$ be a group of finite Morley rank, and let $E \leq G$ be a
finite elementary abelian 2-group of rank at least 3.
Let $\theta$ be a connected nilpotent $E$-signalizer functor.
Then $\theta$ is complete and $\theta(E)$ is nilpotent.
\end{namedtheorem}

We shall work with the proper group $\theta(E)$ in the same manner as
 we did with $\Gamma_E(G)$ in Proposition \ref{GITcondA}.

We also recall a variation on Fact \ref{CgenerationKU2}.

\begin{fact}[{\cite[\qFact 3.7]{Bu03}}]\label{Cgeneration}
Let $H$ be a solvable $p^\perp$-group of finite Morley rank.
Let $E$ be a finite elementary abelian $p$-group acting definably on $H$.
Then  $$ H = \gen{ C_H(E_0) \mid E_0 \leq E, [E:E_0]=p } $$
\end{fact}

\begin{proof}[of Theorem \ref{twogencore}]
It suffices to show that $\Gamma_{E,2}(G) < G$
 by Proposition \ref{twogencore_blowup}.
The Nilpotent Signalizer Functor Theorem
 says that $\theta$ is complete and $\theta(E)$ is nilpotent.
Since $G$ is simple, our result will follow from
$$ \Gamma_{E,2}(G) \leq N_G(\theta(E))\mathperiod $$

For any $U,V\in \elabpgrps_2(E)$,
 we have that $\theta(U) = \theta(V)$ because
$$ \theta(u) \leq \gen{C_{\theta(u)}(v) : v\in V^\#} \leq \theta(V)
  \quad\text{for any $u\in U$\mathperiod} $$
by Fact \ref{Cgeneration} and the signalizer functor property.
Thus $\theta(U) = \theta(E)$.
For any $U\in \elabpgrps_2(E)$ and any $g\in N_G(U)$,
our hypothesis ($\dag$) yields
$$ \theta(E)^g = \theta(U)^g = \theta(U^g) = \theta(U) = \theta(E)\mathperiod $$
Thus $\Gamma_{E,2}(G) \leq N_G(\theta(E)) < G$, as desired.
\end{proof}

\section{Balance and components}\label{sec:reBnC}

In the tame setting of \cite{Bo95},
 Borovik states that $O(C_G(i))$ is a nilpotent signalizer functor.
In view of Theorem \ref{twogencore},
it then follows that either $G$ has a proper 2-generated core, if $O(C_G(i)) \neq 1$,
 or else $C_G(i)$ is ``reductive'' in the sense of Fact \ref{KgrpO}.
It also follows, from Proposition \ref{GITcondA}, that either
 $G$ has a proper 2-generated core, or else $\Gamma_{\Omega_1(S^\o)}(G) = G$.
These two facts constitute the reductivity and generation conditions
 of the Generic Trichotomy Theorem \cite{BB04},
so a tame version of the Generic Trichotomy Theorem then follows.
In fact, \cite{BB04} uses these two conditions merely to establish
that $G$ is generated by the quasisimple components of the centralizers
of toral involutions, and hence by their root $\SL_2$-subgroups.
The remainder of the argument focuses on these root $\SL_2$-subgroups,
treating them as an abstract family of root $\SL_2$-subgroups for $G$,
and eventually applying the Curtis-Tits Theorem.

In this section, we turn our attention towards ``unbalanced groups''
where the group $O(C_G(i))$ is not necessarily a signalizer functor,
in order to eliminate the hypothesis of tameness.  Instead, we use the ``most
unipotent'' parts of $O(C_G(i))$ as signalizer functors.
In Theorem \ref{rebalancing},  these signalizer functors are used
 to prove a dichotomy between a proper 2-generated core, and
 our $\reB$-property (see \S\ref{subsec:rebalancing} below).
Corollary \ref{KgrpO_reU} then provides a limited form of the
 reductivity proved in Fact \ref{KgrpO}.
However, this weaker form of reductivity does not admit such a
quick proof of generation by components.  So our version of this result,
Theorem \ref{rebalanced_Asar2} below, requires a considerably more
delicate argument.

\subsection{Partial balance}\label{subsec:rebalancing}

We require an example to explain the failure of balance.

\begin{example}\label{rebalancing_example}
Consider a field $(k,T,+,\cdot)$ of finite Morley rank,
with $T < k^*$ torsion-free, and $G := \SO_8(k)$ ($D_4$).
By Table 4.3.1 on p.~145 of \cite{GLS3},
there are involutions $i,j$ in $G$, lying in a common torus, such that
$$ C_G(i) \cong \SL_4(k) * k^* \quad\textrm{and}\quad
 C_G(j) \cong \SL_3(k) * \SL_3(k) $$
So $O(C_G(i)) = O(k^*) = T \neq 1$ and $O(C_G(j)) = 1$.
However, every inner involutive automorphism of $\SL_n$
 is a central product with one copy of $k^*$.
So $O(C_G(\cdot))$ is not a signalizer functor.
In fact, the reductivity hypothesis of Fact \ref{KgrpO} fails too,
 although its conclusions still holds since the centralizer is still reductive.
\end{example}

Our solution to this is to choose a reduced rank $\rr^*(\cdot)$
 which is the largest possible problematic reduced rank in $k^*$,
and work above it by using the fact that $\rk(k^*) > \rr_0(k^*)$.

\begin{definition}\label{def:rbarstar}
Consider a simple $K^*$-group $G$ of finite Morley rank and
let $X$ be a subgroup of $G$ with $m(X)\geq3$.
We write $$ I^0(X) := \{ i\in I(X) : m(C_X(i))\geq3 \} $$
for the set of involutions from eight-groups in $\elabpgrps_3(X)$.
We define $$ \rr^O(X) := \sup \{ \rr_0(O(C_G(i))) : i\in I^0(X) \} $$
as the supremum of the reduced ranks of the odd parts of the
centralizers in $G$ of the involutions in $I^0(X)$.
We also define $\rr^*(X)$ to be the supremum of $\rr_0(k^*)$ as $k$ ranges
over the base fields of the quasisimple components of the quotients
$C^\o_G(i) / O(C_G(i))$ associated to involutions $i\in I^0(X)$.
\end{definition}

One can easily check that
$$ \rr^O(G) = \max_{E\in \elabpgrps_3(G)} \rr^O(E)
   \quad\textrm{and}\quad 
   \rr^*(G) = \max_{E\in \elabpgrps_3(G)} \rr^*(E) \mathperiod $$
We recall that, for a nonsolvable group $L$ of finite Morley rank,
$\Uir(L)$ and $U_p(L)$ need not be solvable, as quasisimple algebraic
groups are generated by the unipotent radicals of their Borel subgroups.
Proposition \ref{control_field} below will shed further light on the
definition of $\rr^*(\cdot)$ by providing a converse to this observation.

\begin{definition}\label{def:reE}
We continue in the notation of Definition \ref{def:rbarstar}.
% Let $G$ be a simple $K^*$-group of finite Morley rank, and
%  let $X$ be a subgroup of $G$ with $m(X)\geq3$. 
For a definable subgroup $H$ of $G$, we define $\reU_X(H)$ to be
 the subgroup of $H$ generated by $U_p(H)$ for $p$ prime
 as well as by $\Uir(H)$ for $r > \rr^*(X)$.
As an abbreviation, we use $\reF_X(H)$ to denote $F^\o(\reU_X(H))$,
and $\reE_X(H)$ to denote $E(\reU_X(H))$.
We use $\recomponents^X_Y$ to denote the set of components
of $\reE_X(C_G(i)) = E(\reU_X(C_G(i)))$ for $i\in I^0(Y)$ with $Y \leq X$,
 and we set $\recomponents_X = \recomponents^X_X$.
\end{definition}

$\reU_X(H)$ is the subgroup of $H$ which is generated by its
unmistakably unipotent subgroups.  These definitions are all
sensitive to the choice of $X$, which is usually a fixed eight-group.

\begin{definition}\label{def:reB}
We say that a simple $K^*$-group $G$ with $m(G)\geq3$ satisfies
 the {\em $\reB$-property} if, for every 2-subgroup $X \leq G$ with $m(X)\geq3$
 and every $t\in I^0(X)$, the group $\reU_X(O(C^\o_G(t)))$ is trivial.
This is equivalent to
\begin{conclusions}
\item[($\reB$-1)]  $U_p(O(C_G(t))) = 1$
  for all $t\in I^0(G)$ and every prime $p$.
\item[($\reB$-2)]  $\rr^O(X) \leq \rr^*(X)$
  for every 2-subgroup $X \leq G$ with $m(X)\geq3$, and
\end{conclusions}
\end{definition}

The $\reB$-property is an unbalanced alternative to Borovik's % strong
$B$-conjecture: that $O(C_C(i)) = 1$ for all $i\in I(G)$.
Although the $\reB$-property is significantly more delicate
than the strong $B$-property, the next two subsections will
establish results about the components in $\recomponents_X$
which are similar to Borovik's.

Our goal in this subsection is to verify that the failure of the
 $\reB$-property leads to a proper weak 2-generated core.
For this, we need two appropriate signalizer functors.

\begin{lemma}\label{rebalanced_Uq}
Let $G$ be a simple $K^*$-group of finite Morley rank and odd type
with $m(G)\geq3$, and let $E \in \elabpgrps_3(G)$.
Then $\{ U_p(O(C_G(t))) \mid t \in E^\# \}$
 is a connected nilpotent $E$-signalizer functor satisfying
\begin{texteqn}{\tag{$\dag$}}
$\theta(s)^g = \theta(s^g)$ for all $s\in E^\#$ and
 all $g\in G$ for which $s^g\in E$.
\end{texteqn}
\end{lemma}

We need the following two facts.

\begin{fact}[{\cite[Ex.~11 p.~93 \& Ex.~13c p.~72]{BN}}]\label{nodeftorsion}
Let $G$ be a group of finite Morley rank and let $H \normal G$ be a
definable subgroup.  If $x\in G$ is an element such that $\bar{x}\in G/H$
is a $p$-element, then $x H$ contains a $p$-element.  Furthermore,
if $H$ and $G/H$ are $p^\perp$-groups, then $G$ is a $p^\perp$-group.
\end{fact}

\begin{fact}[{\cite{ABCCdraft}; \cite[\qFact 3.2]{Bu03}}]\label{Cquotient}
Let the group $G = H \rtimes T$ be a semidirect product of finite Morley rank.
Suppose $T$ is a solvable $\pi$-group of bounded exponent and
$Q \normal H$ is a definable solvable $T$-invariant $\pi^\perp$-subgroup.
Then $$ C_H(T)Q/Q = C_{H/Q}(T)\mathperiod $$
\end{fact}

\begin{proof}[of Lemma \ref{rebalanced_Uq}]
Let $\theta(t) := U_p(O(C_G(t)))$.
We observe that $\theta(s)^g = \theta(s^g)$
 for every involution $s\in I(G)$ and every $g\in G$.
By Fact \ref{Upnilpotence}, $\theta(s)$ is nilpotent.
$\theta(s)$ is connected by definition.

Let $s,t\in E^\#$; in particular $[s,t]=1$.
Also let $K_s = O(C_G(s))$.
Since $C_G(s)$ is a $K$-group, Fact \ref{KgrpO} says
$C^\o_G(s)/K_s = G_1 * \cdots * G_n * F$
is the central product of finitely many quasisimple algebraic groups
$G_1,\ldots,G_n$ and of a definable divisible abelian group $F$.
Since $F$ is abelian, $O(C^\o_F(t)) = 1$ by Fact \ref{nodeftorsion}.
We now consider the action of $t$ on the components.
For any component $G_k$, either $t$ normalizes $G_k$, or else
 $t$ swaps $G_k$ with another component $G_l = G_k^t$.
In the second case, the centralizer of $t$ is some diagonal
 subgroup of $G_k^t * G_k$, i.e.\ 
 $C^\o_G(t) \cong \{ (g,\sigma(g)) \mid g\in G_k \} \cong G_k$
 for some automorphism $\sigma$ of $G_k$.
So we may assume that $t$ normalizes each $G_k$ with $k \leq m$,
 and $C^\o_{G_{m+1} * \cdots * G_n}(t) \cong G_{n+m \over 2} * \cdots * G_n$.
By Facts \ref{autalg} and \ref{Creductive2},
 $C^\o_{G_k}(t)$ is reductive for $k \leq m$.
So $O(C^\o_{G_k}(t))$ is a subgroup of an algebraic torus,
 and hence divisible abelian.
Hence $U_p(O(C^\o_{G_1 * \cdots * G_n}(t))) = 1$.
So $U_p(O(C^\o_{C^\o_G(s)/K_s}(t))) = 1$.
Since $C_{C^\o_G(s)}(t) K_s/K_s = C_{C^\o_G(s)/K_s}(t)$
 by Fact \ref{Cquotient},
  $U_p(O(C_{C^\o_G(t)}(s))) = U_p(O(C_{C^\o_G(s)}(t))) \leq K_s$.

For any $t\in E^\#$, the group $\theta(t)$ is $p$-unipotent,
 and so $2^\perp$ and nilpotent.
By Fact \ref{Cconnected}, $C_{\theta(t)}(s)$ is a connected.
So $C_{\theta(t)}(s) = U_p(O( C_{\theta(t)}(s) ))
 \leq U_p(O( C_{C^\o_G(t)}(s) )) \leq U_p(K_s) = \theta(s)$.
\end{proof}

\begin{lemma}\label{rebalanced_Uz}
Let $G$ be a simple $K^*$-group of finite Morley rank and odd type
with $m(G)\geq3$.  Let $E \in \elabpgrps_3(G)$ and set $r := \rr^O(E)$.
If $r > \rr^*(E)$ then
 $\Uir(O(C_G(t)))$ is a connected nilpotent $E$-signalizer functor
 again satisfying $(\dag)$.
\end{lemma}

\begin{proof}%[of Lemma \ref{rebalanced_Uz}]
Let $\theta(t) := \Uir(O(C_G(t)))$.
We observe that $\theta(s)^g = \theta(s^g)$
 for every involution $s\in I(G)$ and every $g\in G$.
$\theta(s)$ is clearly connected and solvable.
So $\theta(t)$ is nilpotent by Theorem \ref{nilpotence}.  % Needs $r := \rr^O(E)$.

Let $s,t\in E^\#$; in particular $[s,t]=1$.
Also let $K_s = O(C_G(s))$.
Since $C_G(s)$ is a $K$-group, Fact \ref{KgrpO} says
$C^\o_G(s)/K_s = G_1 * \cdots * G_n * F$
is the central product of finitely many quasisimple algebraic groups
$G_1,\ldots,G_n$ and of a definable divisible abelian group $F$.
Since $F$ is abelian, $O(C^\o_F(t)) = 1$.
We next show that $\Uir(O(C^\o_{G_1*\cdots*G_n}(t))) = 1$.
For any component $G_k$, either $t$ normalizes $G_k$, or else
 $t$ swaps $G_k$ with another component $G_l = G_k^t$.
In the second case, the centralizer of $t$ is some diagonal
 subgroup of $G_k^t * G_k$, i.e.\ 
 $C^\o_G(t) \cong \{ (g,\sigma(g)) \mid g\in G_k \} \cong G_k$
 for some automorphism $\sigma$ of $G_k$.
So we may assume that $t$ normalizes each $G_k$ with $k \leq m$,
 and $C^\o_{G_{m+1} * \cdots G_n}(t) \cong G_{n+m \over 2} * \cdots G_n$.
Consider a connected definable indecomposable abelian
 subgroup $A$ of $O(C^\o_{G_1 * \cdots * G_m}(t))$
 with $A/J(A)$ torsion-free.
Then there is a $k \leq m$ with
 a nontrivial projection map $\pi : A \to O(C^\o_{G_k/Z(G_k)}(t))$.
By \cite[\qLemma 2.9]{Bu03}, % Fact \ref{Uhom}
 the image $\pi(A)$ is also indecomposable abelian,
 and $\rr_0(A) = \rr_0(\pi(A))$.
By Facts \ref{autalg} and \ref{Creductive2},
 $C^\o_{G_k/Z(G_k)}(t)$ is reductive, and hence
$O(C^\o_{G_k/Z(G_k)}(t))$ is an algebraic torus.
It follows that
$$ \rr_0(A) = \rr_0(\pi(A)) \leq
 \rr_0(O(C^\o_{G_k/Z(G_k)}(t))) \leq \rr^*(E)\mathperiod $$
Since $r > \rr^*(E)$,  we have
 $\Uir(O(C^\o_{G_1*\cdots*G_n}(t))) = 1$.
Since $C_{C^\o_G(s)}(t) K_s/K_s = C_{C^\o_G(s)/K_s}(t)$
 by Fact \ref{Cquotient},
  $\Uir(O(C_{C^\o_G(t)}(s))) = \Uir(O(C_{C^\o_G(s)}(t))) \leq K_s$.

For any $t\in E^\#$, the centralizer $C_{\theta(t)}(s)$ is a connected
$(0,r)$-unipotent $2^\perp$-group by Fact \ref{Cunipotent}. 
Thus $$ C_{\theta(t)}(s) \leq \Uir(O(C_{C^\o_G(t)}(s)))
 \leq \Uir(K_s) = \theta(s)\mathcomma $$
 and $\theta$ is a signalizer functor.
\end{proof}

We can now verify the $\reB$-property,
 in the absence of a proper 2-generated core.

\begin{theorem}\label{rebalancing}
Let $G$ be a simple $K^*$-group of finite Morley rank and odd type
with $m(G) \geq 3$.  Then either
\begin{conclusions}
\item $G$ has a proper weak 2-generated core,  or else
\item $G$ satisfies the $\reB$-property, i.e.\ $\reU_X(O(C_G(t)) = 1$
 for every 2-subgroup $X \leq G$ with $m(X)\geq3$ and every $t\in I^0(X)$.
\end{conclusions}
\end{theorem}

\begin{proof}
We first suppose that ($\reB$-1) fails,
 i.e.\ $U_q(O(C_G(i))) \neq 1$ for some involution $i\in I^0(G)$.
There is an $E\in \elabpgrps_3(G)$ containing $i$.
By Lemma \ref{rebalanced_Uq},
 $\theta(t) := U_p(O(C_G(t)))$ is a connected nilpotent $E$-signalizer functor
 satisfying $(\dag)$.
So $G$ has a proper weak 2-generated core by Theorem \ref{twogencore}.

We next suppose that ($\reB$-2) fails,
 i.e. $\rr^O(X) > \rr^*(X)$ for some 2-subgroup $X \leq G$ with $m(X)\geq3$.
There is an $E\in \elabpgrps_3(X)$ such that $\rr^O(E) = \rr^O(X)$.
Let $\theta(t) := \Uir(O(C_G(t)))$ where $r := \rr^O(E)$ is the
 largest reduced rank appearing inside $O(C(i))$ for involutions $i\in E^\#$.
Now $\theta(t) := U_{0,r}(O(C_G(t)))$ is a connected nilpotent
 $E$-signalizer functor satisfying $(\dag)$ by Lemma \ref{rebalanced_Uz}.
By the choice of $r$, $\theta(i)$ is nontrivial for some involution $i\in E^\#$.
%  and $\theta$ is nilpotent by Fact \ref{nilpotence}.
So $G$ again has a weak proper 2-generated core by Theorem \ref{twogencore}.
\end{proof}

\subsection{Existence of components in $\recomponents_X$}\label{subsec:rebalanced_Asar1}

In this subsection, we will use the $\reB$-property to show that
$G$ is a group ``of component type'' in the sense that 
$$ \hbox{$\recomponents_X \neq \emptyset$
   for every 2-subgroup $X \leq G$ with $m(X)\geq3$.} $$
We employ the $p$-unipotent and 0-unipotent signalizer functors
found in \S\ref{subsec:rebalancing}, via Theorem \ref{rebalancing}.
Our major tool will be the following analog of Fact \ref{KgrpO}
which allows us to exploit the $\reB$-property.

\begin{lemma}\label{KgrpO_U}
Let $H$ be a $K$-group of finite Morley rank and odd type, and
let $\tilde{H}$ be the subgroup of $H$ generated by 
\begin{hypotheses}
\item $\Uir(H)$ for $r > \rr_0(O(H))$, as well as
\item $U_p(H)$ for any prime $p$ satisfying $U_p(O(H)) = 1$.
\end{hypotheses}
Then $\tilde{H} = E(\tilde{H}) * F^\o(\tilde{H})$
 and $F^\o(\tilde{H})$ is abelian.
\end{lemma}

% ...
%
% \begin{fact}[\Altinel and Cherlin \cite{AC99}]\label{centralext}
% A group of finite Morley rank which is a perfect ($G=G'$) central
% extension of a quasisimple algebraic group over an algebraically closed
% field is an algebraic group and has finite center.
% \end{fact}

\begin{proof} % [of Lemma \ref{KgrpO_U}]
We first show that $[F(O(H)),\tilde{H}]=1$.
Let $A$ be a definable connected nilpotent subgroup of $H$
 with either
\begin{enumerate}
\item[a.] $\Uir(A) = A$ for some $r > \rr_0(O(H))$, or
\item[b.] $U_p(A) = A$ for some prime $p$ for which $U_p(O(H)) = 1$.
\end{enumerate}
Then $A \cdot F(O(H))$ is nilpotent
 by Fact \ref{quicknilpotence2}.
Since $r > \rr_0(O(H))$, we have
 $[A, F(O(H))] = 1$ by Fact \ref{nildecomp}.
% As such $A$ generate $\tilde{H},
So $[F(O(H)),\tilde{H}]=1$.

Since $F(O(\tilde{H}))$ is definably characteristic in $\tilde{H}$,
 we have $F(O(\tilde{H})) \leq F(O(H))$, and so
 $F(O(\tilde{H})) \leq Z(\tilde{H})$.
Hence $O(\tilde{H})$ is nilpotent, and
\begin{equation} \tag{$\star$}
 O(\tilde{H}) = F(O(\tilde{H})) \leq Z(\tilde{H}) \mathperiod
\end{equation}

Consider $K := \tilde{H} / O(\tilde{H})$.
By Fact \ref{KgrpO}, $K = E(K) * F(K)$ and $F(K)$ is abelian.
By ($\star$), the inverse image of $F(K)$ in $\tilde{H}$ is nilpotent,
 and thus equals $F^\o(\tilde{H})$.
As any quasi-simple component $L$ of $E(K)$ is perfect,
 such a component $L$ admits only finite central extensions
 by \cite{AC99}.% Fact \ref{centralext}.
By ($\star$), the inverse image $\hat L$ of $L$ in $\tilde{H}$
 is isomorphic to a central product $L * O(\tilde{H})$,
 which thus contains a component of $\tilde{H}$.
Now $\tilde{H} = E(\tilde{H}) * F^\o(\tilde{H})$, as desired.

Therefore $F^\o(\tilde{H})$ satisfies our generation hypotheses.
Clearly $U_p( F^\o(\tilde{H}) ) \le U_p( O(H) ) = 1$.
By Fact \ref{nildecomp},
 $F^\o(\tilde{H})$ is a central product of the $\Uir( F^\o(\tilde{H}) )$
 with $r > \rr_0(O(H)) \geq \rr_0(O(\tilde{H})$.
It follows that $F^\o(\tilde{H})$ is abelian since $F(K)$ was abelian.
\end{proof}
% As $F(K)$ is abelian, $\hat F$ is abelian by Fact \ref{Unilcommutator}.

The $\reB$-property states that the centralizers of appropriate involutions
 satisfy the hypotheses of Lemma \ref{KgrpO_U}.

\begin{corollary}\label{KgrpO_reU}
Let $G$ be a simple $K^*$-group of finite Morley rank and odd type
with $m(G)\geq3$ which satisfies the $\reB$-property.
Then, for every 2-subgroup $X \leq G$ with $m(X)\geq3$ and every $i\in I^0(X)$,
we have $\reU_X(C_G(i)) = \reE_X(C_G(i)) * \reF_X(C_G(i))$
 and $\reF_X(C_G(i))$ is abelian.
\end{corollary}

\smallskip

We can now verify that $\recomponents_X$ is nonempty.

\begin{theorem}\label{rebalanced_Asar1}
Let $G$ be a simple $K^*$-group of finite Morley rank and odd type
with $m(G) \geq 3$.  Then either
\begin{conclusions}
\item $G$ has a proper weak 2-generated core, or else
\item $\recomponents_X \neq \emptyset$
 for every 2-subgroup $X \leq G$ with $m(X)\geq3$.
\end{conclusions}
\end{theorem}

\begin{proof}
By Theorem \ref{rebalancing},
 we may assume that $G$ satisfies the $\reB$-property.
Consider a 2-subgroup $X \leq G$ with $m(X)\geq3$.
There is an $E\in \elabpgrps_3(X)$ with $\rr^*(E)$ maximal.
So $\rr^*(E) = \rr^*(X)$ and $\recomponents_E \subset \recomponents_X$.

We first consider the case where $C^\o_G(i)$ is solvable for all $i\in E^\#$.
In particular, $\rr^*(E) = 0$.
For all $i\in E^\#$,
 $U_p(C_G(i)) = U_p(O(C_G(i)) = 1$ and $\rr^O(E) \leq \rr^*(E) = 0$,
 since $G$ satisfies the $\reB$-property.
By Fact \ref{Udichotomy}, $O(C^\o_G(i))$ is a good torus,
 and hence central in $C^\o_G(i)$ by \cite[Theorem 6.16]{BN}.
Since $C^\o_G(i)/O(C^\o_G(i))$ is abelian by Fact \ref{KgrpO}.
 $C^\o_G(i)$ is nilpotent, and divisible.
Now $C^\o_G(i)'$ is torsion-free by \cite[Theorem 6.9]{BN}.
As $\rr^O(E) = 0$, $C^\o_G(i)$ is in fact abelian.
By Fact \ref{Cgeneration},
 $C^\o_G(i) = \gen{C^\o_{C^\o_G(i)}(E_0) : E_0 \leq E, [E:E_0]=2}$.
As $C^\o_G(i) \neq 1$ by \cite[Ex. 13 \& 15 p.~79]{BN}, % Fact \ref{centinv},
 there is some four-group $E_1 \leq E$ with $H := C^\o_G(E_1) \neq 1$.
Since each $C^\o_G(i)$ is abelian,
 $H \normal C^\o_G(i)$ for all $i\in E_1^\#$,
 and thus $H \normal \Gamma_{E_1}(G)$.
Since $G$ is simple, $\Gamma_{E_1}(G) \leq N_G(H) < G$,
 and $G$ has a proper weak 2-generated core $\Gamma^0_{S,2}(G) < G$
 by Proposition \ref{GITcondA}.
So we may assume that $C^\o_G(i)$ is nonsolvable for some $i\in E^\#$.

We now fix an $i\in E^\#$ and a component $L$ of
 $C^\o_G(i) / O(C^\o_G(i))$ so that $\rr_0(k^*) = \rr^*(E)$
 where $k$ is the base field of $L$.
Since $k$ is algebraically closed, $k^*$ contains torsion, and hence
 $$ \rk(k_+) = \rk(k^*) > \rr_0(k^*) = \rr^*(E)\mathperiod $$
Suppose toward a contradiction that $\reE_E(C_G(i)) = 1$.
By Corollary \ref{KgrpO_reU},  we have
 $\reU_E(C_G(i)) = \reE_E(C_G(i)) * \reF_E(C_G(i))$
 and $\reF_E(C_G(i))$ is abelian for every $i\in E^\#$.
So $\reU_E(C_G(i)) = \reF_E(C_G(i))$ is abelian.
If $\Char(k)>0$ then
 $U_{\Char(k)}(C^\o_G(i)) \leq \reF_E(C^\o_G(i))$ is abelian.
If $\Char(k)=0$ then
 $U_{0,\rk(k_+)}(C^\o_G(i)) \leq \reF_E(C_G(i))$ is abelian.
Either case contradicts the existence of $L$.
\end{proof}

\smallskip

We also observe that the definition of $\recomponents_X$ restricts
 the fields involved as follows.

\begin{proposition}\label{control_field}
Let $H$ be a group of finite Morley rank which is isomorphic to
 a linear algebraic group over an algebraically closed field $k$.
Then
\begin{conclusions}
\item If $\Uir(H) \neq 1$ for some $r > \rr_0(k^*)$
   then $\Char(k) = 0$ and $\rk(k) = r$.
\item If $U_p(H) \neq 1$ then $\Char(k) = p$.
\end{conclusions}
If $H$ is quasisimple, these conditions imply
 $U_p(H) = H$ and $\Uir(H) = H$, respectively.
\end{proposition}

\begin{proof}
If $\Char(k) \neq p$, then $H$ has bounded $p$-rank,
so the second point follows.

We now consider the first point.
Let $A$ be a nontrivial $(0,r)$-unipotent abelian group, and
 let $\hat{A}$ be the Zariski closure of $A$.
Then $\hat{A} = S \times U$ where $S$ is semisimple,
 and $U$ is the unipotent radical of $A$.
If $A$ has nonunipotent elements, then $\bar{A} := A U/U$
 is a nontrivial subgroup of the semisimple group $\hat{A}/U$.
As $\hat{A}/U$ is linear,
 $\hat{A}/U \hookrightarrow (k^*)^n$ for some $n$.
But $\Uir(\bar{A}) = \bar{A}$ by \cite[\qLemma 2.11]{Bu03},
 contradicting $r > \rr_0(k^*)$.
So $A$ consists of unipotent elements, i.e.\ $A \leq U$ ($= \hat{A}$).
Hence $\Char(k) = 0$.
As $U$ is linear, $U \hookrightarrow (k_+)^n$ for some $n$.
By \cite[\qCorollary 3.3]{Po87},
 there are no definable subgroups of $k^+$.
So $\Uir(U) = 1$ unless $r = \rk(k)$.

The last remark follows from the fact that quasisimple algebraic groups
 are generated by the unipotent radicals of their Borel subgroups,
 or indeed by any conjugacy class of elements.
\end{proof}

\subsection{Generation by components in $\recomponents_X$}\label{subsec:rebalanced_Asar2}

We next show that components in $\recomponents_X$ generate $G$, i.e.
$$ \hbox{$\gen{\recomponents_{\Omega_1(S^\o)}} = G$ when $\pr(S) \geq 3$.} $$
However, these results will be proven in a form usable also when $\pr(S) < 3$.

For any group $H$ of finite Morley rank, any 2-subgroup $X$ acting
definably on $H$, and any $V\in \elabpgrps^0_2(X)$, we define
$$ \reGamma_{X,V}(H) = \gen{\reU_X(C^\o_H(v)) : v\in V^\#}\mathperiod $$

\begin{lemma}\label{rebalanced_Asar2_L}
Let $G$ be a $K^*$-group of finite Morley rank and odd type.
Let $X$ be a 2-subgroup of $G$ with $m(X)\geq3$.
Suppose that $G$ satisfies the $\reB$-property and
that there is a four-group $E\in \elabpgrps^0_2(X)$ which
centralizes a Sylow\o 2-subgroup $T$ of $G$.
Let $H := \gen{\reE_X(C_G(z)) : z\in E^\#}$.
Then the following hold.
\begin{conclusions}
\item  For any $x,y\in E^\#$,  we have
  $[\reF_X(C_G(x)),\reF_X(C_G(y))]=1$ and \\
  \hspace*{24pt} the group $\reF_X(C_G(x))$ normalizes $\reE_X(C_G(y))$.
\item $\reU_X(O(H)) = 1$.
\item $\reGamma_{X,E}(G) \leq N^\o_G(\reE_X(H))$. 
\end{conclusions}
\end{lemma}

\begin{proof} %[of Lemma \ref{rebalanced_Asar2_L}]
As a notational convenience,  let $F_x := \reF_X(C_G(x))$ and
 $H_x := \reE_X(C_G(x))$ for  $x\in E^\#$.
Since $G$ satisfies the $\reB$-property, Corollary \ref{KgrpO_reU}
 says $\reU_X(C_G(x)) = H_x * F_x $ and $F_x$ is abelian.
Since $y\in E^\#$ normalizes $F_x$ for any $x\in E^\#$,
 there is a homomorphism $f_y : F_x \to F_x$ given by $u \mapsto [y,u]$.
Since $E$ centralizes $T$ and $F_x \cap T$ is the Sylow 2-subgroup of $F_x$,
 there is no 2-torsion in $F_x \setminus \ker(f_y)$.
So $[F_x,E] \leq O(F_x) \leq O(C^\o_G(x))$
 by Fact \ref{nodeftorsion} (and Fact~\ref{concommutator}).
By Fact \ref{Unilcommutator}, $[\Uir(F_x),E]$ is a $\Uir$-group.
Clearly $[U_p(F_x),E]$ is $p$-unipotent.
Since $\reU_X(C^\o_G(x)) = H_x * F_x$,
 we have $F_x = \reU_X(F_x)$.
So $[F_x,E] \leq \reU_X(O(C^\o_G(x))) = 1$ by the $\reB$-property.
Now $F_x \leq \reU_X(C_G(y))$ for every $x,y\in E^\#$.
Since $F_y$ is central in $\reU_X(C_G(y))$, we have $[F_x,F_y] = 1$.
Also $F_x$ normalizes $H_y$, as $H_y \normal \reU_X(C_G(y))$.
% (see \cite[\qLemma 7.1iii]{BN}).  % was: Fact \ref{components}.
Thus $F := \gen{ F_z : z\in E^\# }$ is an abelian group
 of automorphisms of $H = \gen{ H_z : z\in E^\# }$.

Since $H_x$ is characteristic in $C_G(x)$ for all $x\in E^\#$,
 $E$ normalizes $H$.
Suppose towards a contradiction that $\reU_X(O(H)) \neq 1$.
So either
\begin{enumerate}
\item[(1)] $K := \Uir(O(H)) \neq 1$ for some $r > \rr^*(X)$, or else
\item[(2)] $K := U_p(O(H)) \neq 1$ for some prime $p$.
\end{enumerate}
In either case, $K = \Gamma_E(K)$ by Fact \ref{CgenerationKU2}.
So a $K_x := C_K(x)$ is nontrivial.
In case (1), we may choose $r$ maximal,
 so $K$ is nilpotent by Theorem \ref{nilpotence}.
Hence $K_x$ is $(0,r)$-unipotent by Lemma \ref{Cunipotent}.
In case (2), $K$ is nilpotent by Fact \ref{Upnilpotence}.
Hence $K_x$ is $p$-unipotent by Fact \ref{Cconnected}.
In either case, $K_x$ is nilpotent and normalized by $H_x$.
Since $K_x \leq \reU_X(C_G(x))$, and $H_x$ is semisimple, 
 we have $K_x \leq F_x$, in contradiction to the $\reB$-property.
Thus $\reU_X(O(H)) = 1$, as desired.

For the last part, we may assume $H < G$ is a $K$-group.
So $H = \reE_X(H) * \reF_X(H)$ by Corollary \ref{KgrpO_reU}.
Since $\reE_X(H)$ is characteristic in $H$, $F$ normalizes $\reE_X(H)$.
So $\reGamma_{X,E}(G) = F \cdot H \leq N^\o_G(\reE_X(H))$.
\end{proof}

We now provide the promised variant of Lemma \ref{preGITcondA} above.

\begin{lemma}\label{rebalanced_preGITcondA}
Let $G$ be a $K^*$-group of finite Morley rank and odd type with $m(G)\geq3$,
 and which satisfies the $\reB$-property.
Let $T$ be a Sylow\o 2-subgroup of $G$, and let $X\leq C_G(T)$ be a 2-subgroup
 which centralizes $T$ and has $m(X)\geq3$.
Then $\reGamma_{X,U}(G) = \reGamma_{X,V}(G)$
 for any two four-groups $U,V \in \elabpgrps^0_2(X)$.
\end{lemma}

We need the following algebraic fact.
% version of Theorem \ref{rebalanced_Asar2} below.

\begin{fact}\label{CEgeneration}
Let $G$ be a quasisimple algebraic group over an algebraically closed field
 of characteristic \noteql2 which is not of type $\pPSL_2$, and
let $V$ be a four-group of algebraic automorphisms of $G$
 which centralizes a maximal 2-torus of $G$.
Then $$ G = \gen{ E(C_G(x)) : x\in V^\# }\mathperiod $$
\end{fact}

\begin{proof}
Let $T$ be a maximal 2-torus of $G$ centralized by $V$.
For $i\in V^\#$, Fact \ref{Creductive2} says $C^\o_G(i)$ is reductive.
Hence $C^\o_G(i) = F_i * H_i$ where $H_i := E(C^\o_G(i))$
 and $F_i := F^\o(C^\o_G(i))$ is an algebraic torus.
So $F_i$ is the Zariski closure of $F_i \cap T$,
 and hence is centralized by $V$.
Now $F := \gen{ F_x \mid x\in V^\# }$ is an algebraic torus
 normalizing $H = \gen{ H_x \mid x\in V^\# }$.

We now show that $H$ is reductive.
Suppose that $H$ has a nontrivial unipotent radical $U$.
By Fact \ref{CgenerationKU2},
 $U_j := C_U(j) \neq 1$ for some $j\in V^\#$.
But $U_j \normal C^\o_G(j)$,
 contradicting the reductivity of $C^\o_G(j)$.
So $H$ must be reductive.

Now $\Gamma_V(G) = F H \leq N^\o_G(E(H))$.
We observe that $H_j \neq 1$ for $j\in V^\#$,
 since $G \not\cong \pPSL_2$.
So $H$ is not solvable, and $E(H) \neq 1$.
Since $\Gamma_V(G) = G$ by Fact \ref{CgenerationKU2},
 it follows that $E(H) = G$ because $G$ is quasisimple,
 and hence $H = G$.
\end{proof}

\smallskip

\begin{proof}[of Lemma \ref{rebalanced_preGITcondA}]
We may assume that $[U,V] = 1$
 since $U$ and $V$ lie in the same 2-connected component.
It is enough to show that $\reU_X(C_G(u)) \leq \reGamma_{X,V}(\reU_X(C_G(u))$
 for any $u\in U^\#$.
Since $G$ satisfies the $\reB$-property, % we have
 $\reU_X(C_G(u)) = \reE_X(C_G(u)) * \reF_X(C_G(u))$
 by Corollary \ref{KgrpO_reU}.
By Fact \ref{nildecomp}, the abelian group
 $F_u := \reF_X(C_G(u))$ may be written as a product of various
 $U_p(F_u)$, with $p$ prime, and $\Uir(F_u)$, with $r > \rr^*(X)$.
By Fact \ref{Cunipotent}, $C_{\Uir(F_u)}(v)$ is $(0,r)$-unipotent.
By Fact \ref{Cconnected}, $C_{U_p(F_u)}(v)$ is $p$-unipotent.
So $C_{F_u}(v) \leq \reU_X(C_G(u))$.
Fact \ref{CgenerationKU2} yields
$$ \reF_X(C_G(u)) \leq \reGamma_{X,V}(\reU_X(C_G(u)) \mathperiod $$
Now consider a component $L \normal \reE_X(C_G(u))$.
It suffices to show that $L = \reGamma_{X,V}(L)$.

In the case that $L \not\cong \pPSL_2$, Fact \ref{CEgeneration} yields
$$ L = \gen{ E(C_L(x)) : x\in V^\# }\mathperiod $$
By Proposition \ref{control_field},  we have
 $E(C_L(x)) = \reE_X(C_L(x))$, and hence $L = \reGamma_{X,V}(L)$.

In the case that $L \cong \pPSL_2$, $\pPSL_2$ has no graph automorphisms,
 so every $x\in V^\#$ acts by some inner automorphism by Fact \ref{autalg}.
We observe that $T \leq C_G(u)$, and thus
 $T$ is a Sylow\o 2-subgroup of $C_G(u)$.
So $L \cap T$ is a Sylow\o 2-subgroup of $L$.
Since $\pPSL_2$ contains no four-group centralizing a torus,
 there is now some $x\in V^\#$ which centralizes $L$,
 and the claim follows.
\end{proof}

We now prove a version of Proposition \ref{GITcondA}.

\begin{proposition}\label{rebalanced_GITcondA} 
Let $G$ be an infinite simple $K^*$-group of finite Morley rank and odd type
with $m(G)\geq3$.  Let $T$ be a Sylow\o 2-subgroup of $G$, and
let $X\leq C_G(T)$ be a 2-subgroup which centralizes $T$ with $m(X)\geq3$.
Suppose that $\reGamma_{X,E}(G) < G$ for some $E\in \elabpgrps^0_2(X)$.
Then $G$ has a proper weak 2-generated core.
\end{proposition}

\begin{proof}
By Proposition \ref{GITcondA},
 it is enough to show that $\Gamma_E(G) < G$.
Let $A \in \elabpgrps_3(X)$ be an eight-subgroup of $X$ containing $E$.
By Lemma \ref{weak_CconM}, $\Gamma_E(G) \leq \Gamma_{A,2}(G)$.
So the result will follow from the following claim and simplicity.
$$ \Gamma_{A,2}(G) \leq N_G(\reGamma_{X,E}(G)) $$

By Lemma \ref{rebalanced_preGITcondA} and Fact \ref{Asch46.2}-1,
 $\reGamma_{X,E}(G) = \reGamma_{X,U}(G)$ for any $U \leq A$.
For any four-group $U\leq A$,
 $$ N_G(U) \leq N_G(\reGamma_{X,U}(G)) = N_G(\reGamma_{X,E}(G))\mathperiod $$
Thus $\Gamma_{A,2}(G) \leq N_G(\reGamma_{X,E}(G))$, as desired.
\end{proof}

We now prove that our components generate $G$.

\begin{theorem}\label{rebalanced_Asar2}
Let $G$ be a $K^*$-group of finite Morley rank and odd type with $m(G)\geq3$.
Suppose that there is a four-group $E\in \elabpgrps^0_2(G)$
 which centralizes a Sylow\o 2-subgroup $T$ of $G$, and that
 there is an eight-group $X\in \elabpgrps_3(C_G(T))$ containing $E$.
Then either
\begin{conclusions}
\item $G$ has a proper weak 2-generated core, or
\item $\gen{\recomponents^X_E} = \gen{\reE_X(C_G(x)) : x\in E^\#} = G$.
%  for some $X\in \elabpgrps_3(G)$ which contains $E$.
\end{conclusions}
\end{theorem}

We need the following fact about involutive automorphisms of algebraic groups,
which follows immediately from Table 4.3.1 on p.~145 of \cite{GLS3} and
Fact \ref{autalg}.

\begin{fact}\label{autalg_E}
Let $G$ be a quasisimple algebraic group over an algebraically closed field
of characteristic \noteql2 and
let $\alpha$ be a definable involutive automorphism of $G$.
If $G \not\cong \pPSL_2$, then $E(C_G(\alpha)) \neq 1$. 
\end{fact}

\begin{proof}[of Theorem \ref{rebalanced_Asar2}]
We may assume that $G$ satisfies the $\reB$-property
 by Theorem \ref{rebalancing}.

We first show that $\reE_X(C^\o_G(z)) \neq 1$ for some $z\in E^\#$
 (note $m(E) \geq 2$).
We may assume that $\recomponents_X \neq 1$
 by Theorem \ref{rebalanced_Asar1}, so
there is a component $L \leq \reE_X(C_G(x)) \neq 1$ for some $x\in X^\#$.
We observe that $E$ normalizes $L$
 since $X$ is an eight-group containing $E$.
So any $z\in E^\#$ acts on $L$ via an algebraic automorphism,
 by Fact \ref{autalg}, and hence
 $C^\o_L(z)$ is reductive by Fact \ref{Creductive2}.
By Proposition \ref{control_field},
 $\reU_X(U) = U$ for any unipotent group in $L$.
As quasi-simple groups are generated by their unipotent subgroups,
 $\reU_X(C_G(z)) \geq \reE_X(C_L(z)) = E(C^\o_L(z))$.
If $L \not\cong \pPSL_2$, then
 $E(C^\o_L(z)) \neq 1$ by Fact \ref{autalg_E}, 
and hence $\reE_X(C_G(z)) \neq 1$ by Corollary \ref{KgrpO_reU}.
So we may assume $L \cong \pPSL_2$.
Since $\pPSL_2$ has no graph automorphisms,
 any $z\in E^\#$ acts via inner automorphism, by Fact \ref{autalg}.
Since $\pPSL_2$ contains no four-group centralizing a torus,
 there is now some $z\in E^\#$ which centralizes $L$,
 and $\reU_X(C_G(z))$ is nonabelian.
By Corollary \ref{KgrpO_reU},
 $\reU_X(C_G(z)) = \reE_X(C_G(z)) * \reF_X(C_G(z))$
 and $\reF_X(C_G(z))$ is abelian.
So $\reE_X(C^\o_G(z) \neq 1$ for some $z\in E^\#$.

Let $H := \gen{ \reE_X(C_G(x)) : x\in E^\# }$.
Since $\reU_X(O(H)) = 1$ by Lemma \ref{rebalanced_Asar2_L},
 $H = \reE_X(H) * \reF_X(H)$ and $\reF_X(H)$ is abelian
  by Corollary \ref{KgrpO_reU}.
Since $\reE_X(C^\o_G(x)) \neq 1$, we have $\reE_X(H) \neq 1$.

Since $E$ centralizes $T$, Lemma \ref{rebalanced_Asar2_L} says that
 $\reGamma_{X,E}(G) \leq N^\o_G(\reE_X(H))$.
Therefore $\reGamma_{X,E}(G) \leq N^\o_G(\reE_X(H)) < G$,
 and the theorem follows from Proposition \ref{rebalanced_GITcondA}.
\end{proof}

% \begin{corollary}\label{rebalanced_Asar2_C}
% Let $G$ be a $K^*$-group of finite Morley rank and odd type.
% Suppose that $m(C_S(S^\o))\geq3$ where $S$ is a Sylow 2-subgroup of $G$.
% Then either
% \begin{conclusions}
% \item $G$ has a proper weak 2-generated core $\Gamma^0_{S,2}(G) < G$, or
% \item there is a four-group $A\in \elabpgrps^0_2(G)$ such that
% $$ \gen{\recomponents_A} = G\mathperiod $$
% \end{conclusions}
% \end{corollary}

% % \begin{proof}
% % By Lemma \ref{Ecenttorus},
% %  there is a four-group $A\in \elabpgrps^0_2(G)$
% %  which centralizes a maximal 2-torus $S^\o$.
% % So the result follows from Theorem \ref{rebalanced_Asar2},
% % \end{proof}

\section{The Generic Trichotomy Theorem}\label{sec:generictrico}

We now turn our attention toward proving the following,
 our main result.

\begin{namedtheorem}{Generic Trichotomy Theorem}\label{generictrico}
Let $G$ be a simple $K^*$-group of finite Morley rank and odd type
with $\pr(G) \geq 3$. Then either
\begin{conclusions}
\item $G$ has a proper 2-generated core, i.e.\ $\Gamma_{S,2}(G) < G$, or else
\item $G$ is an algebraic group over an algebraically closed field of characteristic \noteql2.
\end{conclusions}
\end{namedtheorem}

Our strategy is to replicate the proof by Berkman and Borovik of the
 Generic Identification Theorem \cite{BB04},
being careful to use only ``safe'' components,
under the assumption that (1) does not occur.
So we adopt the following standing hypotheses and notation.

\begin{hypothesis}\label{generictrico_hypothesis}
We consider a simple $K^*$-group $G$ of finite Morley rank and odd type
with $\pr(G) \geq 3$, and fix a Sylow 2-subgroup $S$ of $G$.
We also suppose that the 2-generated core of $G$ is not proper,
 i.e.\ $\Gamma_{S,2}(G) = G$.
\end{hypothesis}

% \noindent Only facts shall be exempt from Hypothesis \ref{generictrico_hypothesis}.

As $n(2) \geq \pr(G) \ge 3$,
 we have $\Gamma^0_{S,2}(G) = \Gamma_{S,2}(G) = G$ by Fact \ref{Asch46.2}-2.
So, by Theorem \ref{rebalancing},
% the condition $\Gamma^0_{S,2}(G) = G$ implies that
\begin{texteqn}{}
 $G$ satisfies the $\reB$-property.
\end{texteqn}

\subsection{Root $\mathrm{SL}_2$subgroups}\label{subsec:RootGroups}

The first stage in our analysis is to select, and establish the properties of
 a family of abstract ``root $\SL_2$-subgroups'' of $G$.
% To understand the role of these root $\SL_2$-subgroups, the reader may
% wish to glance ahead at conditions (L1) through (L6) of Fact \ref{Lyons}.
%
The root $\SL_2$-subgroups of an {\em algebraic group} associated to
a maximal torus $T$ may be defined as those {\em Zariski closed}
subgroups of $G$ which are normalized by $T$ and are isomorphic
to $\pPSL_2$, or alternatively in terms of groups generated by
opposite root groups.
We employ several facts about root $\SL_2$-subgroups of algebraic groups.

\begin{fact}[{cf. \cite[\qFact 2.1]{BB04}}]\label{rootSL2pairs}
Let $G$ be a quasisimple algebraic group over an algebraically closed field.
Let $T$ be a maximal torus in $G$ and let $K,L$ be {\em Zariski closed}
subgroups of $G$ that are isomorphic to $\SL_2$ or $\PSL_2$ and are
normalized by $T$.  Then
\begin{conclusions}
\item Either $K$ and $L$ commute or
  $\gen{K,L}$ is a quasisimple algebraic group of type $A_2$, $C_2$, or $G_2$.
\item The subgroups $K$ and $L$ are root $\SL_2$-subgroups of $\gen{K,L}$.
\item If $\gen{K,L}$ is of type $G_2$, then $G = \gen{K,L}$.
\end{conclusions}
\end{fact}

More generally, a semisimple subgroup of a simple algebraic group $G$
 which is normalized by a maximal torus $T$ is called
 {\em subsystem} subgroup of $G$, associated to $T$.
Berkman and Borovik refer to the full classification of
 semisimple subsystem subgroups \cite[2.5]{seitz83}
 (see also \S3.1 of \cite{seitz}) for the proof of this fact.
The elementary argument here is based on the following fact.

\begin{fact}[{\cite[\qProposition 3.1]{seitz}; \cite[2.5]{seitz83}}]
\label{subsystem}
Let $G$ be a simple algebraic group, let $T$ be a maximal torus of $G$, and
let $X$ be a closed connected subgroup of $G$ which contains $T$.
Then $X = D Z U$ where $D$ is a subsystem subgroup normalized by $T$,
 $Z$ is a torus, and $U$ is the unipotent radical of $X$.
\end{fact}

\begin{proof}[of Fact \ref{rootSL2pairs}]
By Fact \ref{subsystem}, $\gen{K,L} = D Z U$ where
 $D$ is a subsystem subgroup, $Z$ is a torus, and
 $U$ is the unipotent radical of $\gen{K,L}$.
There is an automorphism $\phi$ of the root system
 for $G$ which sends any root $\alpha \in I$ to
 its negative $-\alpha$, and
 $\phi$ translates to an automorphism $\phi$ of
 the group $G$ such that $\phi$ normalizes $T$
 and $\phi(X_\alpha) = X_{-\alpha}$ \cite{Carter89}.
Since $K$ and $L$ each contain one positive and one negative root from $I$,
we find that $K$, $L$, and $\gen{K,L}$ are all normalized by $\phi$.
If $U$ is nontrivial, it contains a root group $X_\alpha$.
Since $U$ is characteristic in $\gen{K,L}$,
 the nilpotent group $U$ must also contain $X_{-\alpha}$,
 and hence contains a copy of $\pPSL_2$ \cite{Carter89}.
So $U=1$, and $\gen{K,L} = D$.  Since $D$ is semisimple,
$D \cong A_1 * A_1$, $A_2$, $C_2$, or $G_2$, as desired.
\end{proof}

\begin{fact}[{see \cite[p.~19]{Carter93}}]\label{rootSL2}
Let $G$ be a semisimple algebraic group over an algebraically closed field,
and fix a maximal algebraic torus $T$ of $G$.  Then the following hold.
\begin{conclusions}
\item $G$ is generated by its root $\SL_2$-subgroups associated with $T$.
\item The intersection $T \cap K$ of $T$ and a root $\SL_2$-subgroup $K$
 associated to $T$ is a maximal algebraic torus of $K$.
% \item The Weyl group $W := N_G(T)/T$ of $G$ is generated by the Weyl groups
%   $N_{T L}(T)/T$ of the root $\SL_2$-subgroups associated to $T$.
\end{conclusions}
\end{fact}

\begin{proof}
The fact that $G$ is generated by those root $\SL_2$-subgroups
which are normalized by $T$ can be found on p.~19 of \cite{Carter93}.
For the second part, we observe that the maximal algebraic tori of $N_G(K)$,
one of which is $T$ and one of which extends a maximal algebraic torus
of $K$, are conjugate in $N_G(K)$.
\end{proof}

We also need to know that a root $\SL_2$-subgroup is ``cut out'' by the
 centralizer of a 2-torus in the associated maximal torus,
We remark that this is an essential point if one hopes to apply
 Fact \ref{rootSL2pairs}, but it remains somewhat obscure in \cite{BB04}.

\begin{fact}\label{rootSL2cutout}
Let $G$ be a quasisimple algebraic group over
 an algebraically closed field of characteristic \noteql2.
Let $T$ be a maximal algebraic torus of $G$, and let $L$ be
 a root $\SL_2$-subgroup of $G$ normalized by $T$.
Then $L = E(C_G( C_{S^\o}(L) ))$ where $S^\o$ is the Sylow\o 2-subgroup of $T$.
\end{fact}

\begin{proof}
We may assume that $G$ is not isomorphic to $\pPSL_2$
 because otherwise $G = L$.
% Let $L$ be a root $\SL_2$-subgroup of $G$ associated to $T$.
Let $S$ be a Sylow 2-subgroup of $G$ such that $S^\o \leq T$.
Any connected definable group of automorphisms of $G$
 must be inner by Fact \ref{autalg}.
Since $L$ is normalized by $S^\o$,
 we have $\pr(C_{S^\o}(L)) = \pr(G) - 1$.
By Fact \ref{Creductive2} (see also \cite[\qTheorem 3.5.4]{Carter93}),
 $C_G(C_{S^\o}(L))$ is reductive.
So $\pr(K) \leq 1$ where $K := E(C_G(C_{S^\o}(L)))$.
Since $L \leq K$ and $L$ and $K$ are both algebraic subgroups,
 we have $L = K$.
\end{proof}

% We will also need the following model theoretic result, due to Poizat.
%
% \begin{fact}[]\label{PoizatSL2}
% Let $G$ be a group of finite Morley rank which is isomorphic to $\pPSL_2$
% over an algebraically closed field.
% Then $G$ is the only definable connected nonsolvable subgroup of $G$.
% \end{fact}

\smallskip

We now proceed with the analysis of groups satisfying
 Hypothesis \ref{generictrico_hypothesis}.

\begin{lemma}\label{generictrico_maxalgtori}
For any $i\in \Omega_1(S^\o)^\#$ and any definable connected quasisimple % semisimple
 algebraic $L \leq E(C^\o_G(i))$ which is normalized by $S^\o$, we have
\begin{conclusions}
% \item $L$ is normalized by $S^\o$.
\item $S^\o \cap L$ is a Sylow\o 2-subgroup of $L$,
\item $T_L := C_L(S^\o \cap L)$ is a maximal algebraic torus of $L$,
\item $S^\o = C^\o_{S^\o}(L) (S^\o \cap L)$,
   and $\pr(G) = \pr(S^\o) = \pr(C^\o_{S^\o}(L)) + \pr(S^\o \cap L)$.
\end{conclusions}
\end{lemma}

For this, we need the following fact about algebraic groups.

\begin{fact}\label{Sylow_algmaxtori}
Let $G$ be a quasisimple algebraic group over an algebraically closed field
of characteristic \noteql2, and let $T$ be a Sylow\o 2-subgroup of $G$.
Then $C_G(T)$ is a maximal algebraic torus of $G$.
\end{fact}

\begin{proof}
In an algebraic group over an algebraically closed field,
 the maximal algebraic torus is the Zariski closure of $T$,
 and is thus centralized by anything centralizing $T$.
But a maximal algebraic torus is self-centralizing
 by \cite[24.1]{Hump}.
So the result follows.
\end{proof}

\begin{proof}[of Lemma \ref{generictrico_maxalgtori}]
Since $S^\o$ is a Sylow\o 2-subgroup of $N^\o_G(L)$,
 the group $S^\o \cap L$ is a Sylow\o 2-subgroup of $L$.
By Fact \ref{Sylow_algmaxtori},
 $T_L$ is a maximal algebraic torus of $L$.
By Fact \ref{autalg},
 the connected definable group $d(S^\o)$ acts by inner automorphisms on $L$,
 so the third condition follows.
\end{proof}

\begin{lemma}\label{generictrico_recomponents}
For any component $M\in \recomponents_{S^\o}$,
 $M$ is normalized by $S^\o$.
\end{lemma}

\begin{proof}
Clearly $M \normal C^\o_G(i)$ for some $i\in \Omega_1(S^\o)$
 (see \cite[\qLemma 7.1iii]{BN}).  % was: Fact \ref{components}.
Since $S^\o \leq C^\o_G(i)$, $M$ is normalized by $S^\o$.
\end{proof}

\begin{definition}\label{def:generictrico_Sigma}
Let $\Sigma$ be the set of all root $\SL_2$-subgroups of components
$K\in \recomponents_{S^\o}$ which are associated to $T_K$,
 i.e.\ $\Sigma$ is the set of all {\em Zariski closed} subgroups of
 the components $K\in \recomponents_{S^\o}$ which are normalized
 by $T_K$ and are isomorphic to $\pPSL_2$.
\end{definition}

Since the root $\SL_2$-subgroups of $K$ generate $K$ by Fact \ref{rootSL2}-1,
Theorem \ref{rebalanced_Asar2} yields the following.
\begin{equation} \tag{($\star$)}
 \gen{\Sigma} = \gen{\recomponents_{S^\o}} = G
\end{equation}

We view the subgroups in $\Sigma$ as abstract root $\SL_2$-subgroups for $G$.

\begin{lemma}\label{generictrico_Sigma1}
For any $L\in \Sigma$, we have
\begin{conclusions}
\item $L$ is normalized by $S^\o$,
\item $L = E(C_G( C_{S^\o}(L) ))$, and
\item $L$ is a Zariski closed subgroup of any definable quasisimple $K< G$
 which contains $L$ and which is normalized by $S^\o$.
\end{conclusions}
\end{lemma}

\begin{proof}
Let $R_L := C_{S^\o}(L)$ and
let $M \in \recomponents_{S^\o}$ be a component containing $L$
 as a root $\SL_2$-subgroup associated with $T_M$.
By Lemma \ref{generictrico_maxalgtori}-3,
 $S^\o = C_{S^\o}(M) (M \cap S^\o)$.
Since $L$ is normalized by $M \cap S^\o \leq T_M$,
 $L$ is normalized by $S^\o$.
By Fact \ref{rootSL2cutout}, $L = E(C_M( M \cap R_L ))$.
Fix $i\in \Omega_1(S^\o)^\#$ with $M \normal C^\o_G(i)$.
Clearly $E(C_G( R_L )) \leq M$ because any other component
 of $C^\o_G(i)$ meets $R_L$ in an infinite 2-torus.
As $i \in R_L$, $L = E(C_G( R_L ))$.

Now for any definable quasisimple $K < G$
 which contains $L$ and which is normalized by $S^\o$,
the group $R_L$ acts on $K$ by inner automorphisms by Fact \ref{autalg},
 so $L = E(C_K(R_L))$ is Zariski closed.
\end{proof}

\begin{lemma}[{cf.\ \cite[\qLemma 3.1]{BB04}}]\label{generictrico_Sigma2}
For any distinct $K,L \in \Sigma$,
\begin{conclusions}
\item $C_{S^\o}(K) \cap C_{S^\o}(L) \neq 1$
  and $M := \gen{K,L}$ is a $K$-group.
\item Either $K$ and $L$ commute or $M$ is an algebraic group
  of type $A_2$, $B_2 = C_2$, or $G_2$.
\item $S^\o \cap M = (S^\o \cap K) * (S^\o \cap L)$
  is a Sylow\o 2-subgroup of $M$.
\item $K$ and $L$ are root $\SL_2$-subgroups of $M$ normalized by $T_M$.
\item $[T_K,T_L] = 1$.
\end{conclusions}
\end{lemma}

\begin{proof}
Let $R_L := C^\o_{S^\o}(L)$.
Since $S^\o$ normalizes $K$ and $L$ by Lemma \ref{generictrico_Sigma1}-1,
 and $K,L \cong \SL_2$, we know that
 $\pr(R_K),\pr(R_L) = \pr(G) - 1$ and $S^\o = R_K R_L$
 by Lemma \ref{generictrico_maxalgtori}-3.
So
\begin{eqnarray*}
\pr(G) &=& \pr(R_K) + \pr(R_L) - \pr(R_K \cap R_L) \\
 &=& 2 \pr(G) - 2 - \pr(R_K \cap R_L)
\end{eqnarray*}
and $$ \pr(R_K \cap R_L) = \pr(G) - 2 \geq 1. $$
Thus $C_{S^\o}(K) \cap C_{S^\o}(L) \neq 1$ and $M$ has a nontrivial center.
Since $G$ is simple, $M < G$ is a $K$-group.

Let $i \in I(C_{S^\o}(K) \cap C_{S^\o}(L))$.
By Corollary \ref{KgrpO_reU} (and the $\reB$-property),
 $K,L \leq \reE_{S^\o}(C_G(i))$.
If the groups belong to different components of $C^\o_G(i)$,
 then they commute.
If they both belong to the same component $H\in \recomponents_{S^\o}$,
 then $H$ is a quasisimple algebraic group normalized by $S^\o$
 by Lemma \ref{generictrico_recomponents}.
By Lemma \ref{generictrico_Sigma1}-3, $K$ and $L$ are Zariski closed in $H$.
By Fact \ref{rootSL2pairs}-1,
 $M = \gen{K,L}$ is an algebraic group of type $A_2$, $B_2 = C_2$, or $G_2$.
In either case, $S^\o$ normalizes $M$, so
 $S^\o \cap M$ is a Sylow\o 2-subgroup of $M$ and
 $T_M$ is a maximal ``algebraic'' torus of $M$
 by Lemma \ref{generictrico_maxalgtori}.

For (4) and (5), we may assume that $[L,K]\neq1$
 and $M$ is a quasisimple algebraic group.
By Fact \ref{rootSL2pairs}-2,
 $K$ and $L$ are root $\SL_2$-subgroups of $M$.
By Fact \ref{rootSL2}-2, $T_M = T_K * T_L$,
 so $[T_K,T_L] = 1$.
\end{proof}

We give $\Sigma$ a graph structure by placing an edge between
 $L,K\in \Sigma$ when $[L,K] \neq 1$.
Since $G$ is simple and $\gen{\Sigma} = G$ ($\star$),
 the graph $\Sigma$ is connected.
By Lemma \ref{generictrico_Sigma2}-2,
 any adjacent $L,K\in \Sigma$ are algebraic groups over
 the same algebraically closed field.
So all the elements of $\Sigma$ are algebraic groups over
 a common algebraically closed field $\F$.
Since $G$ has odd type, $\Char(\F) \neq 2$.
In particular, $\rk(K) = \rk(L)$ for all $K,L\in \Sigma$.

From this point on, our argument reduces to that given by
 Berkman and Borovik in \cite{BB04},
following the presentation of \cite{BBBC07}.
Indeed, we may now lighten our standing hypotheses to the following.

\begin{hypothesis}\label{generictrico_hypothesis_light}
We consider a simple group $G$ of finite Morley rank and odd type
with $\pr(G) \geq 3$, and fix a Sylow 2-subgroup $S$ of $G$.
Also choose {\em some} family $\recomponents$ of algebraic components
 from the centralizers of involutions in $S^\o$.
Let $\Sigma$ be the set of root $\SL_2$-subgroups,
 from components in $\recomponents$, which are associated to $S^\o$,
 in the sense of Definition \ref{def:generictrico_Sigma}.
Suppose that Lemmas \ref{generictrico_maxalgtori},
 \ref{generictrico_Sigma1}, and \ref{generictrico_Sigma2} are satisfied,
and also that
\begin{equation} \tag{($\star$)}
 \gen{\Sigma} = G\mathperiod
\end{equation}
\end{hypothesis}

\noindent We give the analysis in full below.

\subsection{Weyl group}

We now turn our attention to the Weyl group of $G$,
 continuing under Hypothesis \ref{generictrico_hypothesis_light}.

\begin{lemma}\label{generictrico_torus}
The {\em natural} torus $T := \gen{T_L : L\in \Sigma}$
is divisible abelian.  So $T_L = T \cap L$.
\end{lemma}

\begin{proof}
By Lemma \ref{generictrico_Sigma2}-5,
 the algebraic tori $T_K$ for $K\in \Sigma$ all commute,
 so the result follows.
\end{proof}

\begin{definition}
For any $L\in \Sigma$, $W(L) := N_L(T_L)/T_L$
 is the Weyl group of $L$ and has order 2.
We may identify $W(L) \cong N_L(T)/C_L(T)$ with its image
 in $W := N_G(T)/C_G(T)$, by Lemma \ref{generictrico_torus}.
Now let $r_L$ denote the single involution inside $W(L)$,
 and define $W_0 := \gen{ r_L \in W : L\in \Sigma }$.
\end{definition}

\begin{lemma}[{cf.\ \cite[\qLemma 3.5]{BB04}}]\label{generictrico_rKvsK}
For any $L,K\in \Sigma$,
 $[K,L]=1$ if and only if $[r_K,r_L]=1$.
\end{lemma}

\begin{proof}
It suffices to check this in $\gen{K,L}$.
So the result follows from Fact \ref{rootSL2pairs}-2 and Fact \ref{rootSL2}.
\end{proof}

We will analyze $W_0$ by examining its action on $S^\o$ and $T$.

\begin{lemma}[{cf.\ \cite[Lemmas 3.6 \& 3.7]{BB04}}]\label{generictrico_SconT}
$S^\o$ is the Sylow 2-subgroup of $T$ and $C_G(S^\o) = C_G(T)$.
In particular, $W_0$ acts faithfully on $S^\o$.
\end{lemma}

\begin{proof}
By Lemma \ref{generictrico_torus}, $T$ is divisible abelian,
 so its Sylow 2-subgroup is connected by \cite[Theorem 9.29]{BN}. % Fact \ref{Sylow_con}.

Let $D := S^\o \cap T$.  Suppose toward a contradiction that $D < S^\o$.
For all $K\in \Sigma$, $[S^\o,r_K] = S^\o \cap K$,
 so $r_K$ acts trivially on $S^\o/D$.
Let $b\in S^\o$ satisfy $\abs{b/D} \geq 4$ and
 let $a\in S^\o$ satisfy $a^{\abs{W_0}} = b$.
Then $c = \Pi_{w\in W_0} a^w$ satisfies $b/D = c/D$ and $\abs{c}\geq 4$.
Since $S^\o = C_{S^\o}(K) * (S^\o \cap K)$ and
 $\abs{C_{S^\o}(K) \cap (S^\o \cap K)} \leq \abs{Z(K)} = 2$,
 we have $[C_{S^\o}(r_K) : C_{S^\o}(K)] \leq 2$.
Since $c\in C_{S^\o}(r_K)$,
 $c^2 \in C_{S^\o}(K)$ for all $K\in \Sigma$, and $c^2 \neq 1$.
So $c^2 \in C_G(\gen{\Sigma}) = Z(G)$,
 contradicting the simplicity of $G$.
Thus $S^\o$ is the Sylow 2-subgroup of $T$, and $C_G(S^\o) \geq C_G(T)$.

For the reverse direction, consider $x\in C_G(S^\o)$.
Then, for every $L\in \Sigma$,  $x$ centralizes $C_{S^\o}(L)$.
So $x$ normalizes $L = E(C_G(C_{S^\o }(L)))$
 by Lemma \ref{generictrico_Sigma1}-2. % was: and Fact \ref{components}.
Since $x$ centralizes the maximal 2-torus $S^\o \cap L$,
 $x$ must act on $L$ as an element of $T_L$ by Fact \ref{autalg}.
Thus $x\in C_G(T)$ and $C_G(S^\o) \leq C_G(T)$.
\end{proof}

We use the action of $W_0$ on $S^\o$ to obtain a complex representation.

\begin{lemma}[{cf.\ \cite[\S3.3]{BB04}}]\label{generictrico_reflections_Tate}
$W_0$ has a faithful irreducible representation $R$ over $\C$
 of dimension $\pr(G) \geq 3$
in which the $r_L$ act as reflections for $L\in \Sigma$.
\end{lemma}

For this, we employ a Tate module over the 2-adics.

\begin{fact}[{\cite[\S3.3]{Be01,BB04}}]\label{endptor}
Let $T$ be a $p$-torus of \Prufer $p$-rank $n$ in a group of finite Morley rank.
Then $\End(T)$ can be faithfully represented as the ring of $n \times n$
matrices over the $p$-adic integers $\Z_p \cong \End(\Z(p^\infty))$.
\end{fact}

\begin{proof}[of Lemma \ref{generictrico_reflections_Tate}]
For every $L \in \Sigma$ and every 2-torus $X \leq S^\o$ disjoint from $L$,
 $X$ must act on $L$ as elements of $S^\o \cap L$
 by Fact \ref{autalg},
so there is some 2-torus $Y \leq C_{S^\o}(L)$ with $X \leq Y (S^\o \cap L)$.
So $S^\o = C_{S^\o}(r_L) * (S^\o \cap L)$.
By Lemmas \ref{generictrico_Sigma1}-1 and \ref{generictrico_maxalgtori},
 $r_L$ inverts $S^\o \cap L = [S^\o,r_L]$.
Thus $r_L$ acts as a ``reflection'' on $S^\o$.

By Lemma \ref{generictrico_SconT}, $W_0$ acts faithfully on $S^\o$.
By Fact \ref{endptor}, $W_0$ has a faithful representation over
 the ring of 2-adic integers $\Z_2$ which has dimension $\pr(S^\o) \geq 3$,
By tensoring with $\C$, $W_0$ has a faithful representation $R$ over $\C$
 which has dimension $\pr(S^\o) \geq 3$.
The $r_L$'s continue to act as reflections in this representation.

Now suppose towards a contradiction that $W_0$ acts reducibly on $R$.
Since the representation $R$ is completely reducible, $R =  R_1 \oplus R_2$
 where $R_1$ and $R_2$ are proper $W_0$-invariant subspaces.

Suppose that $W_0$ acts trivially on $R_i$.
Then there is a 2-torus $\hat{R}_i$ centralized by all $r_L$,
 and $G \leq C_G(\hat{R}_i)$, a contradiction.
So we may assume that $W_0$ acts non-trivially on both $R_1$ and $R_2$.

For $L\in \Sigma$, the $-1$-eigenspace $[R,r_L]$ of $r_L$ belongs to
one of the two subspaces, either $R_1$ or $R_2$.
So $r_L$ acts as a reflection on that subspace and centralizes the other.
Let $\Sigma_i := \{ L \in \Sigma : [R,r_L] \leq R_i \}$ for $i=1,2$.
For $L\in \Sigma_1$ and $K\in \Sigma_2$,  we have
 $[r_L,r_K] = 1$, and thus $[L,K] = 1$ by Lemma \ref{generictrico_rKvsK},
 in contradiction with the fact that $\Sigma$ is connected.
\end{proof}

To further constrain $W_0$,
 we next obtain representations of $W_0$ over almost all finite fields.

\begin{lemma}[{cf.\ \cite[\S3.4]{BB04}}]\label{generictrico_reflections_q}
For primes $q > \max(\abs{W_0},3)$ with $q \neq \Char(L)$
for any $L\in \Sigma$,  $W_0$ has a faithful irreducible representation
 over $\Z/q\Z$, where, for any $L \in \Sigma$,
 the involution $r_L$ acts by reflection.
\end{lemma}

\begin{proof}
Consider the elementary abelian $q$-group $E_q$
 generated by all elements of order $q$ in $T$.
$W_0$ clearly acts on $E_q$.  Let $N = N_G(T)$.
Since $C_G(T) \leq C_N(E_q)$, we may show that $W_0$ acts faithfully
 by showing that $C_N(E_q) \leq C_G(T)$.

For any $x\in C_N(E_q) \leq N$, $x$ acts on $\Sigma$ by conjugation.
For any $L\in \Sigma$, if $L^x \neq L$ then $L$ and $L^x$ either
 commute or generate a quasisimple group as root $\SL_2$-subgroups
 by Lemma \ref{generictrico_Sigma2}-2.
In either case, $\abs{L \cap L^x} \leq 2$,
 in contradiction to the fact that $L \cap E_q = L^x \cap E_q$.
So $x$ normalizes $L$, and the element
 $x$ acts on $T \cap L$ as an element of $N_L(T \cap L)$ by Fact \ref{autalg}.
Since the Weyl group of $\SL_2$ inverts the torus,
 any element of $N_L(T \cap L) \setminus C_L(T \cap L)$ inverts some element of $E_q$.
So $x$ centralizes $T \cap L$ for all $L \in \Sigma$.
Now $x$ centralizes $T = \gen{ T \cap L | L \in \Sigma }$,
 and $W_0$ acts faithfully on $E_q$.

We also observe that $W_0$ acts by reflections on $E_q$ because,
 for every $L\in \Sigma$, $[E_q,r_L]$ has order $q$ and is inverted by $r_L$,
 i.e.\ $\abs{E_q \cap L} = q$.

Now suppose toward a contradiction that $W_0$ acts reducibly on $E_q$.
Since $q > \abs{W_0}$, the representation is completely reducible,
 and $E_q =  R_1 \oplus R_2$
 where $R_1$ and $R_2$ are proper $W_0$-invariant subspaces of $E_q$.

Suppose that $W_0$ acts trivially on $R_i$.
For any $L\in \Sigma$,
 $R_i$ acts by inner automorphisms on $L$ by Fact \ref{autalg}.
We recall that nontrivial Weyl group elements in $W(L)$ invert the torus $T_L$.
Since $R_i$ centralizes $T_L$, we know that $R_i$ acts via conjugation
 by elements of $T_L$.
Since $W_0$ centralizes $R_i$, we find that $R_i$ centralizes $L$.
So $R_i \leq Z(\gen{\Sigma}) = Z(G)$, a contradiction.
So we may assume that $W_0$ acts nontrivially on both $R_1$ and $R_2$.

For $L\in \Sigma$, the eigenspace $[E_q,r_L]$ of $r_L$ belongs to
one of the two subspaces, either $R_1$ or $R_2$.
So $r_L$ acts as a reflection on that subspace and centralizes the other.
Let $\Sigma_i := \{ L \in \Sigma : [E_q,r_L] \leq R_i \}$ for $i=1,2$.
Then $\Sigma = \Sigma_1 \cup \Sigma_2$.
For $L\in \Sigma_1$ and $K\in \Sigma_2$,  we have
 $[r_L,r_K] = 1$, and thus $[L,K] = 1$ by Lemma \ref{generictrico_rKvsK},
 in contradiction with the fact that $\Sigma$ is connected.
\end{proof}

The two preceding lemmas provide sufficient information
 to identify the Weyl group $W_0$.

\begin{lemma}[{cf.\ \cite[\qLemma 3.11]{BB04}}]\label{generictrico_reflections}
There exists an irreducible root system $I$ of
 type $A_n$, $B_n$, $C_n$, $D_n$, $E_6$, $E_7$, $E_8$, or $F_4$
 on which $W_0$ acts as a crystallographic reflection group.
\end{lemma}

This lemma follows from the following major fact,
 which depends on a detailed analysis of
 the irreducible complex reflection groups \cite{ShTo54,Cohen76}.

\begin{fact}[{\cite[\qTheorem 2.3]{BBBC07}}]\label{reflection:coxeter}
Let $W$ be a finite group, $I\subseteq W$ a subset, and $n$ an
integer, satisfying the following conditions.
\begin{enumerate}
\item The set $I$ generates $W$, consists of involutions, and is
closed under conjugation in $W$; \item The graph $\Delta_I$ with
vertices $I$ and edges $(i,j)$ for noncommuting pairs $i,j\in I$
is connected; \item For all sufficiently large prime numbers
$\ell$, $W$ has a faithful representation $V_\ell$ over the finite
field ${\mathbb{F}}_\ell$ in which the elements of $I$ operate as
complex reflections, with no common fixed vectors.
\end{enumerate}
Then one of the following occurs.
\begin{itemize}
\item[$(a)$] $W$ is a dihedral group  acting in dimension $n=2$,
or cyclic of order two. \item[$(b)$] $W$ is isomorphic to an
irreducible crystallographic Coxeter group, that is,
$A_n,B_n,C_n,D_n$ ($n\ge 3$), $E_n$ ($n=6,7$,or $8$), or $F_n$
($n=4$), \item[$(c)$] $W$ is a semidirect product of a quaternion
group of
  order $8$ with the symmetric group ${\rm Sym}_3$, acting naturally,
represented in dimension $2$.
\end{itemize}
If, in addition, over some field, $W$ has an irreducible
representation of dimension at least $3$, in which the elements of
$I$ act as reflections, then case $(b)$ applies.
\end{fact}

\begin{proof}[of Lemma \ref{generictrico_reflections}]
We observe that $\{ r_L : L \in \Sigma \}$ is a normal subset of $W_0$
 which generates $W_0$.
The noncommuting graph on this set is connected by Lemma \ref{generictrico_rKvsK}.
So Lemmas \ref{generictrico_reflections_q} and \ref{generictrico_reflections_Tate}
 complete the verification of the hypotheses of Fact \ref{reflection:coxeter}.
\end{proof}

\smallskip

We also show that all reflections in $W_0$
 come from our root $\SL_2$-subgroups. 

\begin{lemma}[{cf.\ \cite[\qLemma 3.12]{BB04}}]\label{generictrico_roots}
Every $r\in W_0$ which is a reflection in the representation $R$ over $\C$
 has the form $r_K$ for some $K\in \Sigma$.
\end{lemma}

Recall that the reflections of a Coxeter group correspond to roots in the
associated root system (see \cite[Lemma 5.7]{HumpCoxeter}),
 and hence there are at most two conjugacy classes of reflections.

\begin{fact}[{\cite[10.4 \qLemma C]{HumpLieAlg}}]\label{f:reflection_conjclasses}
A finite irreducible reflection group of type $A_n$, $D_n$, $E_6$, $E_7$, or
$E_8$ has only one conjugacy class of reflections.  A finite irreducible
reflection group of type $B_n$, $C_n$, $F_4$, and $G_2$ has two conjugacy
classes of reflections, corresponding to the short and long roots.
\end{fact}

Since the roots of only one length are closed under the action of the
Coxeter group, they form the root system for a proper subgroup.

\begin{fact}\label{f:reflection_oneconjclass}
The subgroup of $B_n$, $C_n$, $F_4$, or $G_2$ generated by the reflections
associated to roots of only one length is a proper subgroup.
\end{fact}

\begin{proof}[of Lemma \ref{generictrico_roots}]
By Fact \ref{f:reflection_conjclasses},
 there are at most two conjugacy classes of reflections in $I$,
 corresponding to the short and long roots.
So we may assume that $I$ has more than one root length,
 i.e.\ $W_0 \cong B_n$, $C_n$, or $F_4$, and
that the set $S := \{ r_L : L\in \Sigma \}$ consists of
 only one of these conjugacy classes.
By Fact \ref{f:reflection_oneconjclass}, $\gen{S} < W_0$,
 a contradiction.
\end{proof}

\subsection{Identification}

We continue the analysis of the preceding subsections, loosely
following \cite[\S3.6]{BB04}.  We will invoke the
Curtis-Tits theorem which may be expressed as follows:
 a simply connected quasisimple algebraic group is the free amalgam
 of the system of subgroups and inclusion maps corresponding to
 all root $\SL_2$ subgroups and subgroups generated by pairs of
 such subgroups, taken relative to a fixed maximal torus \cite{GLS2}.
The Generic Identification Theorem of Berkman and Borovik proceeds
by passing from the full system of groups and subgroups to the
collection of subsystems corresponding to pairs of roots, which
are now known.   A flexible form of this result is based on a
result of Timmesfeld \cite{Timm04}.
% The following form of this result is based on Timmesfeld's work \cite{Timm04}.

\begin{fact}[{\cite[\qProposition 2.3]{BBBC07}}]\label{CurtisTits}
Let $\Phi$ be an irreducible root system (of spherical type) and
rank at least $3$, and
let $\Pi$ be a system of fundamental roots for $\Phi$.
Let $X$ a group generated by subgroups $X_r$ for $r\in  \Pi$,
Set $X_{rs}=\gen{X_r,X_s}$.
Suppose that $X_{rs}$ is a group of Lie type $\Phi_{rs}$
over an infinite field, with $X_r$ and $X_s$ corresponding
root $\SL_2$-subgroups with respect to some maximal torus of $X_{rs}$.
Then $X/Z(X)$ is isomorphic to a group of Lie type via a map
carrying the subgroups $X_r$ to root $\SL_2$-subgroups.
\end{fact}

% Note that if $X$ is, in addition, a group of finite Morley rank,
%  then $X$ is itself a Chevalley group by Fact \ref{centralext}.

We now conclude the proof of the Generic Trichotomy Theorem,
 working, as usual, under Hypothesis \ref{generictrico_hypothesis_light}.
By Lemma \ref{generictrico_reflections},  $I$ is the desired
 irreducible root system of spherical type and rank at least 3.
For every vertex $i\in I$, there is an $r_i\in W_0$
 which is a reflection in the representation $R$ over $\C$.
There is a also reflection $r_L$ for every $L \in \Sigma$,
By Lemma \ref{generictrico_roots},
 there is an $L_i \in \Sigma$ such that $r_i = r_{L_i}$, and
$$ \gen{ L_i | i\in I } = \gen{\Sigma} = G \mathperiod $$
For $i,j\in I$, the group $M := \gen{L_i,L_j}$ is of Lie type
 by Lemma \ref{generictrico_Sigma2}-2 when $[L_i,L_j] \neq 1$.
If $[L_i,L_j] = 1$ then $M = L_i * L_j$ which has Lie type
 because $L_i$ and $L_j$ are algebraic over the same field.
By Lemma \ref{generictrico_Sigma2}-4,
 $L_i$ and $L_j$ are root $\SL_2$-subgroups corresponding
 to a maximal torus $T_M$ of $M$.
Now $G$ is a Chevalley group by Fact \ref{CurtisTits}, as desired. % \qed

\section*{Acknowledgments}

I thank my adviser Gregory Cherlin for direction during my thesis work,
of which this article is an outgrowth, and for guidance during its later
developments.  I am grateful to Tuna \Altinel for his careful reading of
this material, and numerous corrections and suggestions.  
I am also grateful to my collaborators on \cite{BBN04} and \cite{BCJ},
without whom this article would merely provide the Generic Trichotomy
Theorem, as well as Alexandre Borovik and \Ayse Berkman for advice,
and allowing me to use their arguments from \cite{BB04}.

Financial support for this work comes from
 NSF grant DMS-0100794 and DFG grant Te 242/3-1.
I also thank IGD at Universit\'e Lyon I for their hospitality.

\small
\bibliographystyle{alpha}
\bibliography{burdges,fMr}

% \affiliationone{Jeffrey Burdges\\
% School of Mathematics, The University of Manchester\\
% PO Box 88, Sackville St., Manchester M60 1QD, England\\
% \email{Jeffrey.Burdges@manchester.ac.uk}}

\end{document}